\documentclass{article}
\usepackage{arxiv}
\usepackage{graphicx}
\usepackage{amsmath}
\usepackage{amsfonts}
\usepackage{amssymb}
\usepackage{bm}

\usepackage{subcaption}
\usepackage[normalem]{ulem}
\usepackage{xcolor}

\usepackage[utf8]{inputenc}
\usepackage[T1]{fontenc}    
\usepackage{hyperref}       
\usepackage{url}            
\usepackage{booktabs} 
\usepackage{nicefrac}   
\usepackage{microtype} 
\usepackage{stackengine}
\stackMath

\def\sss{\scriptscriptstyle}
\setstackgap{L}{8pt}

\newcommand{\bx}{\bm{x}}
\newcommand{\dd}{\mathrm{d}}
\newcommand{\As}{\mathsf{A}}
\newcommand{\fs}{\mathsf{f}}
\newcommand{\hs}{\mathsf{h}}
\newcommand{\Ks}{\mathsf{K}}
\newcommand{\Ws}{\mathsf{W}}
\newcommand{\Us}{\mathsf{U}}
\newcommand{\Vs}{\mathsf{V}}
\newcommand{\SGs}{\mathsf{\Sigma}}
\newcommand{\Ss}{\mathsf{S}}
\newcommand{\Po}{\mathrm{P}_{\mathrm{\Sigma}}}

\newcommand{\px}{\mathrm{p}_{\mathfrak{D}}}
\newcommand{\Nx}{{\mathrm{N}_{\mathrm{x}}}}
\newcommand{\Nm}{{\mathrm{N}_{\mathrm{modes}}}}
\newcommand{\Ns}{{\mathrm{N}_{\mathrm{spls}}}}
\newcommand{\lc}{{\ell_c}}

\newcommand\freefootnote[1]{
  \let\thefootnote\relax
  \footnotetext{#1}
}

\title{Numerical Considerations for the Construction of Karhunen-Lo\`{e}ve Expansions}

\author{
 Cosmin Safta \\
  Sandia National Laboratories\\
  Livermore, CA 94550 \\
  \texttt{csafta@sandia.gov} \\
   \And
 Habib N. Najm \\
  Sandia National Laboratories\\
  Livermore, CA 94550 \\
  \texttt{hnnajm@sandia.gov} \\
}

\begin{document}
\maketitle

\begin{abstract}
This report examines numerical aspects of constructing Karhunen-Lo\`{e}ve expansions (KLEs) for second-order stochastic processes. The KLE relies on the spectral decomposition of the covariance operator via the Fredholm integral equation of the second kind, which is then discretized on a computational grid, leading to an eigendecomposition task. We derive the algebraic equivalence between this Fredholm-based eigensolution and the singular value decomposition of the weight-scaled sample matrix, yielding consistent solutions for both model-based and data-driven KLE construction. Analytical eigensolutions for exponential and squared-exponential covariance kernels serve as reference benchmarks to assess numerical consistency and accuracy in 1D settings. The convergence of SVD-based eigenvalue estimates and of the empirical distributions of the KL coefficients to their theoretical $\mathcal{N}(0,1)$ target are characterized as a function of sample count. Higher-dimensional configurations include a two-dimensional irregular domain discretized by unstructured triangular meshes with two refinement levels, and a three-dimensional toroidal domain whose non-simply-connected topology motivates a comparison between Euclidean and shortest interior path distances between the grid points. The numerical results highlight the interplay between the discretization strategy, quadrature rule, and sample count, and their impact on the KLE results.
\end{abstract}


\section{Introduction}
\label{sec:intro}
Random field representations are a foundational tool in computational science and engineering, providing a rigorous mathematical framework for quantifying field uncertainty in physical systems. Applications span a broad range of disciplines, including subsurface flows, material science, and uncertainty quantification (UQ) for high-fidelity simulations. A common challenge in all of these settings is the need to represent, discretize, and sample from infinite-dimensional stochastic processes using a finite number of parameters, while preserving the statistical structure of the underlying field. The Karhunen-Lo\`{e}ve expansion (KLE)~\cite{Karhunen:1946} addresses this challenge by providing an optimal, mean-square-convergent series spectral representation of a second-order stochastic process in terms of a set of deterministic orthogonal basis functions and uncorrelated scalar random variables.\freefootnote{Paper published as Sandia National Laboratories Technical Report SAND2026-18707.}

In its continuous form, the KLE is computed via the solution of a Fredholm integral equation of the second kind, whose kernel is the covariance function of the process. The resulting eigen-pairs - eigenvalues and eigenfunctions - provide a hierarchical decomposition of the process variance, with the leading modes capturing the dominant spatial correlations. The KLE is optimal in the mean-square sense, i.e. any truncation to a finite number of terms minimizes the mean-square reconstruction error among all linear orthogonal expansions of the same dimension, which lends itself useful for dimensionality reduction. This property is shared with principal component analysis (PCA) in the data-driven setting~\cite{Rasmussen:2005}.

The use of KLEs in computational science and UQ gained momentum following the work of Ghanem and Spanos \cite{Ghanem:1991}, who employed it in the stochastic finite element method to propagate input uncertainties through partial differential equation (PDE) models. Their framework, combined with polynomial chaos expansions for the output quantities of interest, established a spectral approach to stochastic computation. Le Ma\^{i}tre and Knio~\cite{OlmOmk:2010} provide a comprehensive overview of spectral methods for UQ that encompasses KLE methods. The accuracy and efficiency of the downstream UQ computation depends on the fidelity of the KLE used to parameterize uncertain input fields.

Numerous studies have examined the numerical aspects of constructing KLEs in practical settings. The discretization of the Fredholm equation on a finite spatial mesh converts the continuous eigenproblem into a discrete one. Betz~{\em et.al.}~\cite{Betz:2014} provide a systematic comparison of discretization strategies - including the Galerkin, collocation, and midpoint methods - and their accuracy as a function of the mesh resolution relative to the correlation length. Schwab and Todor~\cite{Schwab:2006} analyzed the convergence properties of KLE approximations and developed fast multipole-based methods to reduce the computational cost of assembling and solving large covariance eigenproblems. Alternative computational approaches avoid the explicit formation of the covariance matrix by working with an ensemble of field realizations: the singular value decomposition (SVD) of the sample matrix yields the same eigenpairs as the discretized Fredholm equation when quadrature weights are incorporated appropriately. This connection to PCA and the SVD is well-established in the statistical literature~\cite{Rasmussen:2005} and motivates data-driven KLE constructions in settings where an analytical covariance model is not available. The choice of covariance kernel - exponential, squared-exponential, and others - governs the smoothness of the random field realizations and the rate of eigenvalue decay, which in turn determines how many KLE terms are needed for a given accuracy. Analytical eigensolutions exist only for a limited class of kernels on simple domains, e.g. the exponential kernel on a bounded interval~\cite{Ghanem:1991} and the squared-exponential kernel on the real line under a Gaussian spatial measure~\cite{Rasmussen:2005,Zhu:1998}, and these serve as benchmarks for validating numerical implementations.

The accuracy of a numerical KLE depends on the interplay between several modeling choices: the discretization of the computational domain, the correlation length of the covariance kernel relative to the domain size, and - in the data-driven setting - the number of available random field samples. Under-resolved grids fail to capture the oscillatory structure of high-frequency eigenmodes, while an insufficient sample count introduces statistical variability in the estimated eigenvalues and may lead to empirical KLE coefficients that deviate from the expected distribution. These issues become particularly pronounced in multi-dimensional settings and on complex, non-simply-connected domains, where the choice of distance metric used to evaluate the covariance function (Euclidean versus shortest interior path) can affect the eigensolution.

In this technical report we present the theoretical underpinnings of the KLE, derive analytical solutions for the exponential and squared-exponential covariance kernels, and then study the effect of modeling choices highlighted above on the numerical KLE results. The report is organized as follows. Section~\ref{sec:KLEtheory} establishes the mathematical framework, presenting the Karhunen-Lo\`{e}ve expansion, the Fredholm integral equation of the second kind, and the relationship between the continuous expansion and its discretized counterparts, including both the assembly-based and SVD-based solution strategies. Section~\ref{sec:numerical} presents a numerical study of KLE construction across a range of mesh discretizations, correlation lengths, and sample sizes. One-, two-, and three-dimensional test cases are considered, including domains with unstructured meshes and a toroidal geometry for which the effect of the choice of distance metric (Euclidean versus shortest interior path) on the eigensolution is examined. Section~\ref{sec:conclusions} summarizes this report.

\section{Karhunen-Lo\`{e}ve Expansions}
\label{sec:KLEtheory}

\subsection{Definitions}
Let $(\Omega, \Sigma, \Po)$ be a probability space, with $\Omega$ being the sample space, $\Sigma$ a $\sigma-$algebra over $\Omega$, and $\Po$ a probability measure assigning probabilities for each event $\omega$ in $\Sigma$. We consider a $d-$dimensional space $\mathfrak{D}\subset\mathbb{R}^d$.

Let $Y(\bx,\omega)$, with $\bx\in\mathfrak{D}$, e.g. spatial and/or temporal coordinates, and $\omega\in\Sigma$ a random event, be a stochastic process over the space $\mathfrak{D}\times \Omega$, $Y:\mathfrak{D}\times \Omega\rightarrow \mathbb{R}$. For any fixed $\bm{x}\in\mathfrak{D}$, $Y(\bm{x},\cdot)$ is a random variable and for any fixed event $\omega\in\Sigma$, $Y(\cdot,\omega)$ is a sample from a stochastic process over $\mathfrak{D}$. 

We consider the case with second-order random variables (RVs), i.e. RVs with finite second-order moments, equipped with the inner product $\langle\cdot,\cdot\rangle$ over the event space
\begin{equation}
\langle Y_1(\bm{x},\omega),Y_2(\bm{x},\omega)\rangle_{\omega}
=\int_{\Omega}Y_1(\bm{x},\omega)Y_2(\bm{x},\omega)\dd\Po(\omega)=\mathbb{E}_{\omega}[Y_1Y_2].
\end{equation}
Without loss of generality and to simplify notations below, we will assume a zero-mean stochastic process, i.e. one has already subtracted the mean:
\[
Y(\bm{x},\omega)\rightarrow Y(\bm{x},\omega)-\mathbb{E}_{\omega}[Y]
\]
where $\mathbb{E}_{\omega}[Y]=\int_{\Omega}Y(\bm{x},\omega)\dd\Po(\omega)$.  

Let $C(\bx_1,\bx_2)$ denote the covariance function of the stochastic process $Y$. By definition, for an $L_2$ space, the covariance function is symmetric, bounded, and positive definite. Furthermore, by Mercer's theorem, we have that $C$ admits a spectral decomposition in terms of the eigenvalues and eigenfunctions of its corresponding Hilbert-Schmidt integral operator $T_C$
\[
(T_Cf)(\bm{x})=\int_\mathfrak{D} C(\bm{x},\bm{x'}) f(\bm{x'})\dd\mu_\mathfrak{D}(\bx'). 
\]
where $\mu_\mathfrak{D}$ is a Lebesgue measure over $\mathfrak{D}$. Following Ref~\cite{Rasmussen:2005}, here we consider cases where there is a probability density $\px(\bx)$ over $\mathfrak{D}$ and $\dd\mu_\mathfrak{D}(\bx)=\px(\bx)\dd\bx$. 

The eigenvalues and eigenfunctions of the integral operator $T_C$ are the solutions of the Fredholm equation of the second kind
\begin{equation}
\int_{\mathfrak{D}}C(\bx_1,\bx_2) f(\bx_2)\px(\bx_2)\dd\bx_2=\lambda f(\bx_1).
\label{eq:fred}
\end{equation}
The eigenvalues $\lambda_k$ are real, and arranged in decreasing order, $\lambda_1\geq\lambda_2\geq\ldots$. Further, the sum of their squares is finite
\[
\sum_{k=1}^{\infty}\lambda_k^2<\infty.
\]
and
\begin{equation}
C(\bx_1,\bx_2)=\sum_{k=1}^{\infty}\lambda_k f_k(\bx_1) f_k(\bx_2).
\end{equation}

The eigenfunctions $f$ are orthonormal over $\mathfrak{D}$ under the measure $\mu_\mathfrak{D}$,
\begin{equation}
\int_{\mathfrak{D}} f_k(\bx)f_l(\bx)\px(\bm{x})\dd\mu_\mathfrak{D}(\bx)=\delta_{kl}.
\label{eq:kleorth}
\end{equation}
The Karhunen-Lo\`{e}ve expansion of the zero-mean stochastic $Y(\bx,\omega)$ is given by
\begin{equation}
Y(\bx,\omega)=\sum_{k=1}^{\infty}\sqrt{\lambda_k}f_k(\bx)\xi_k(\omega)
=\sum_{k=1}^{\infty}f_k(\bm{x})\zeta_k(\omega).
\end{equation}
In the first sum above $\xi_k$ are uncorrelated RVs with zero-means and unit variances. In the second sum, $\sqrt{\lambda_k}$ is folded inside the RVs. When using this form, $\zeta_k$ are uncorrelated RVs with zero means and variances $\lambda_k$. These RVs are given by
\begin{equation}
\zeta_k(\omega)=\sqrt{\lambda_k} \xi_k(\omega)=\langle Y(\bm{x},\omega) f_k(\bx)\rangle_{\bm{x}}.
\label{eq:rvest}
\end{equation}

\subsection{Analytical Solutions of the Fredholm Equation}
\label{sec:analytical}
In this section we present two cases with covariance kernels that exhibit analytical solutions for the Fredholm equation eigenproblem in Eq.~\eqref{eq:fred}. In the first case, we consider an exponential covariance kernel on a finite 1D domain. In the second case we present results for a squared-exponential covariance kernel on an infinite 1D domain.

\subsubsection{Karhunen–Lo\`{e}ve Expansion for Exponential Covariance Kernel}
\label{sec:klexpanl}

Consider the one–dimensional, stationary exponential covariance kernel
\begin{equation}
C(x_1,x_2) = \exp(-|x_1 - x_2|/\lc),
\label{eq:expcov}
\end{equation}
defined on the finite interval, $x_1,x_2 \in [0,L]$. Here, $\lc$ is the correlation length of the kernel. We seek eigenvalues $\lambda_k$ and orthonormal eigenfunctions $f_k(x)$ satisfying Eq.~\ref{eq:fred} with $\px(x)=1/L$. To simplify the analysis, we will work with normalized coordinates
\[
x_1\rightarrow x_1/L,\,x_2\rightarrow x_2/L,\,\lc\rightarrow \lc/L
\]
using the normalized coordinates and keeping the same notation for the normalized coordinates
\[
\int_0^1 C(x_1,x_2)\,f_k(x_2)\,\dd x_2 = \lambda_k\,f_k(x_1)
\]
This equation admits an analytical solution for the corresponding eigenfunctions and eigenvalues
\[
f_k(x) =
\frac{\sqrt{2}\lc\omega_k}{\sqrt{\lc^2\omega_k^2 + 2\lc+1}}
\left[\cos(\omega_k x) + \frac{1}{\lc\omega_k}\sin(\omega_k x)\right].
\]
and 
\[
\lambda_k = \frac{2\lc}{\lc^2\omega_k^2+1},
\]
where $\omega_k$ is the solution of the equation
\[
\tan\omega_k = \frac{2\lc\omega_k}{\lc^2\omega_k^2-1}.
\]
For more details on the derivation to obtain the above expressions see Ref.~\cite{Ghanem:1991}. 

\begin{figure}[!htb]
\begin{center}
\includegraphics[width=0.45\linewidth]{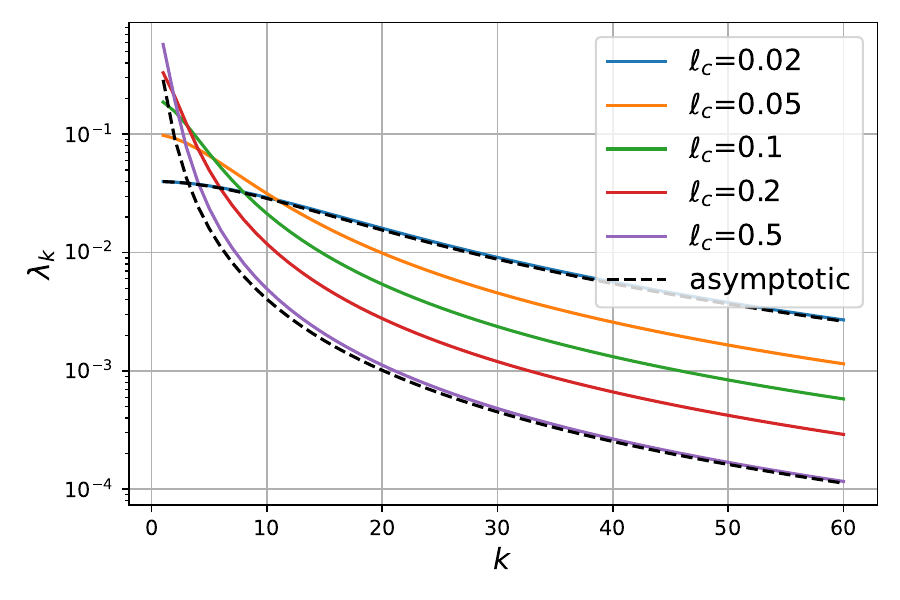}  
\end{center}
\caption{\label{fig:exp1D_evals} Eigenvalues corresponding to a 1D KL expansion with an exponential covariance kernel and a domain size $L=1$. The two dashed lines correspond to asymptotic trends evaluated via Eq.~\eqref{eq:e1Dasymp} for $\lc=0.02$ and $\lc=0.5$, respectively.}
\end{figure}

For large $k$, $\omega_k \approx k\pi$ and 
\begin{equation}
\lambda_k \approx \frac{2\lc}{(k\pi \lc)^2+1}
\approx \left(\frac{2}{\pi^2\lc}\right)\frac{1}{k^2}
\label{eq:e1Dasymp}
\end{equation}
and the eigenfunctions approach pure cosines
\[
f_k(x)\rightarrow \sqrt{2}\cos(k\pi x).
\]
For short correlation lengths, $\omega_k\rightarrow (2k-1)\pi/2$. Figure~\ref{fig:exp1D_evals} shows the analytical eigenvalue spectra for several correlation lengths. Also shown, with dashed lines, are the asymptotic results corresponding to a small correlation length ($\lc=0.05$) and a large correlation length ($\lc=0.5$), relative to the domain size, $L=1$. For the small correlation length, the analytical and asymptotic results are close, verifying the estimates above. Similarly, the asymptotic estimates for large $k$ are verified through the agreement observed for $\lc=0.5$ for $k>30$.

\begin{figure}[!htb]
\begin{center}
\includegraphics[width=0.9\linewidth]{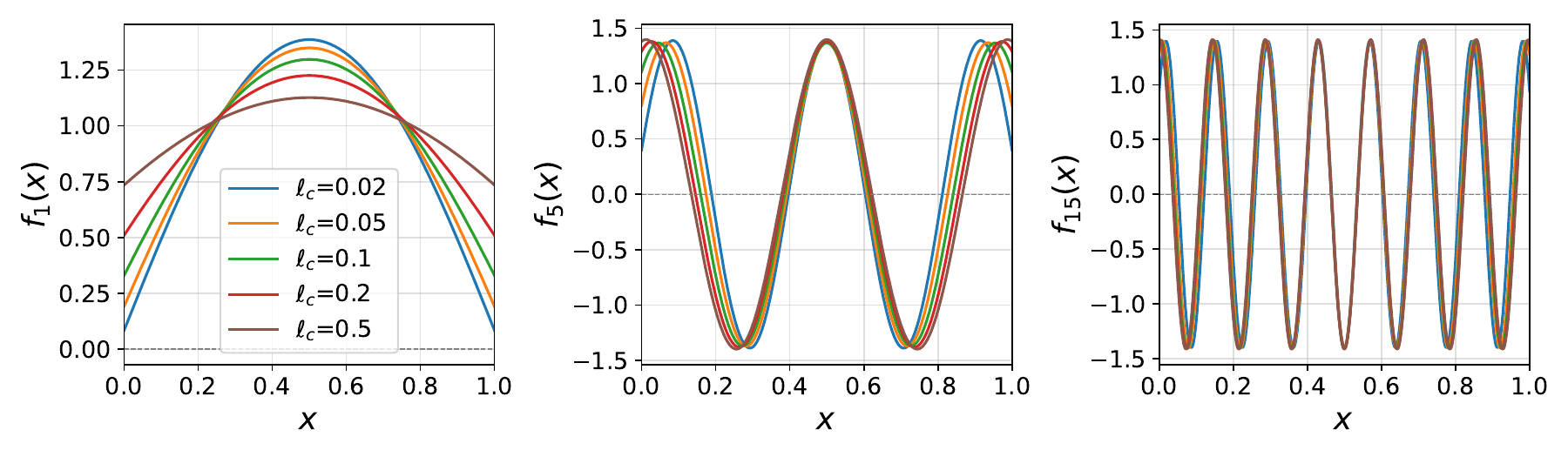}  
\end{center}
\caption{\label{fig:exp1D_evecs} Eigenvectors $f_1$, $f_5$, and $f_{15}$ corresponding to a 1D KL expansion with an exponential covariance kernel and a domain size $L=1$.}
\end{figure}

Figure~\ref{fig:exp1D_evecs} shows eigenfunctions for select eigenmodes and the same choices of correlation lengths as in Figure~\ref{fig:exp1D_evals}. For the first mode, in the left frame, there is a strong dependence on the correlation length, as a direct consequence of small frequency. 
For the higher modes (or frequencies), the eigenfunctions approach the $\sqrt{2}\cos(k\pi x)$ asymptote functional form. 

\subsubsection{Karhunen–Lo\`{e}ve Expansion for Squared-Exponential Covariance Kernel}
\label{sec:anlsqexp}
In this section we explore random fields generated by the covariance kernel
\begin{equation}
C(x_1,x_2) = \exp(-(x_1 - x_2)^2/(2\lc^2)),
\label{eq:sqexpcov}
\end{equation}
over a one-dimensional space, with $x\sim \mathcal{N}(0,\sigma_x^2)$. This case admits analytical results~\cite{Rasmussen:2005,Zhu:1998} for the KL expansion eigenvalues $\lambda_k$ and eigenfunctions $f_k$, as follows. 

\begin{figure}[!htb]
\begin{center}
\includegraphics[width=0.45\linewidth]{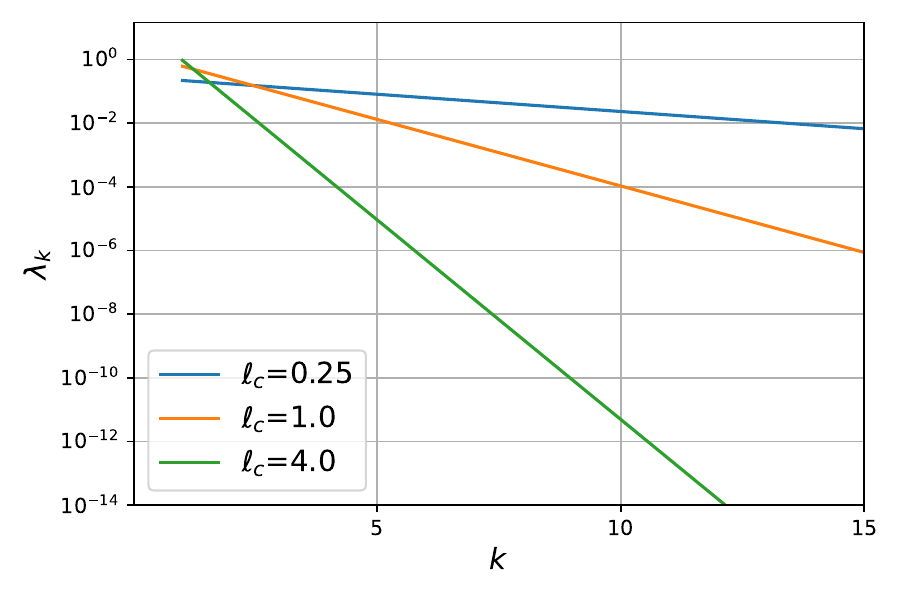}  
\includegraphics[width=0.45\linewidth]{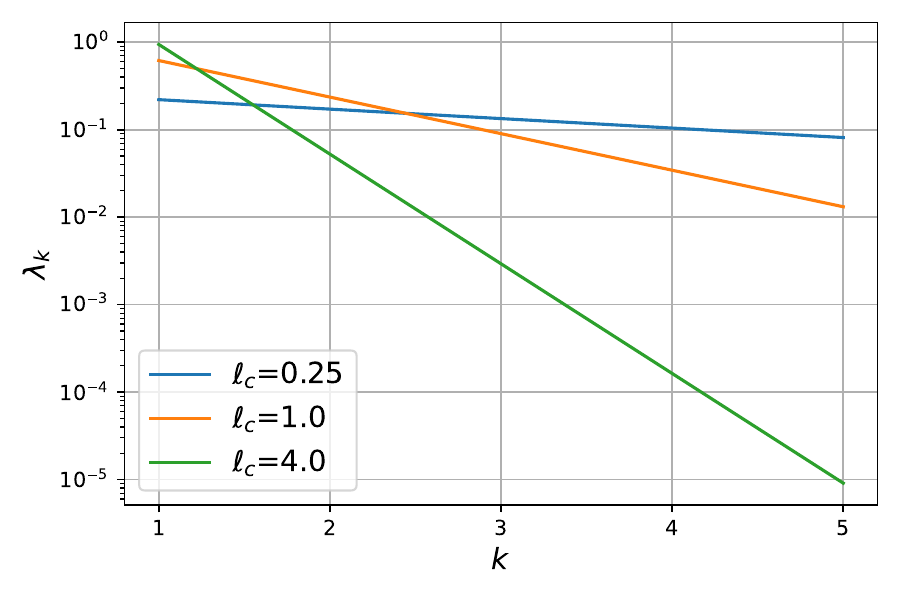}  
\end{center}
\caption{\label{fig:sqexp1D_evals} Eigenvalues of a 1D KLE with a squared-exponential covariance kernel and $x\sim\mathcal{N}(0,\sigma_x)$ with $\sigma_x=1$. Note that correlation lengths are selected as $\lc=\{\sigma_x/4,\sigma_x,4\sigma_x\}$. The right frame shows a detail of the eigenvalue spectra up to the 5-th mode.}
\end{figure}

We will use the same notation as in Ref.~\cite{Rasmussen:2005} 
\[
a=\frac{1}{4\sigma_x^2},\,\,\, b=\frac{1}{2\lc^2},\,\,\,c=\sqrt{a^2+2ab}
\]
and 
\[
A=a+b+c,\,\,\, B=b/A.
\]
The eigenvalues and eigenfunctions are given by
\begin{equation}
\lambda_k = \sqrt{\frac{2a}{A}} B^{k-1},\,\,\, f_k(x)=\sqrt{\frac{2\sigma_x\sqrt{c}}{2^{k-1}(k-1)!}}\exp(-(c-a)x^2)H_{k-1}(\sqrt{2c}x),\,\,\, k=1,2,\ldots
\label{eq:KLgauss_anl}
\end{equation}
where $H_k(x)$ is the $k$-th order ``physicist'' Hermite polynomial, orthogonal with respect to $\exp(-x^2)$
\begin{equation}
\frac{1}{\sqrt{2\pi}k!}\int_{-\infty}^{\infty}H_k(x)H_l(x)\exp(-x^2)dx=\delta_{kl}
\end{equation}
with $H_0(x)=1$, $H_1(x)=2x$, and
\[
H_{k+1}(x) = 2x H_k(x)-2k H_{k-1}(x).
\]
For a given choice of $\sigma_x$ and $\lc$, the expression for the eigenvalues in Eq.~\eqref{eq:KLgauss_anl} depends only on the ratio $\rho=\lc/\sigma_x$. In fact, after some algebraic manipulation,
\begin{equation}
\lambda_k = \frac{2^{k-\frac{1}{2}}\rho}{\left(\rho^2+2+\rho\sqrt{\rho^2+4}\right)^{k-1/2}}.
\label{eq:eval_anl_expression}
\end{equation}
For $\rho=\lc/\sigma_x=1/4$, the eigenvalues decay about $17\%$ per mode and for $\rho=4$ the decay rate is about $94\%$ per mode, with two modes necessary to capture $99\%$ of the process variance.  For the case with $\lc=\sigma_x$, i.e. $\rho=1$,  the decay rate is about $60\%$ per mode. Figure~\ref{fig:sqexp1D_evals} shows the eigenvalue spectra which are straight lines in the semi-log plot given the power law dependence shown in Eq.~\eqref{eq:eval_anl_expression}.

\begin{figure}[!htb]
\begin{center}
\includegraphics[width=0.32\linewidth]{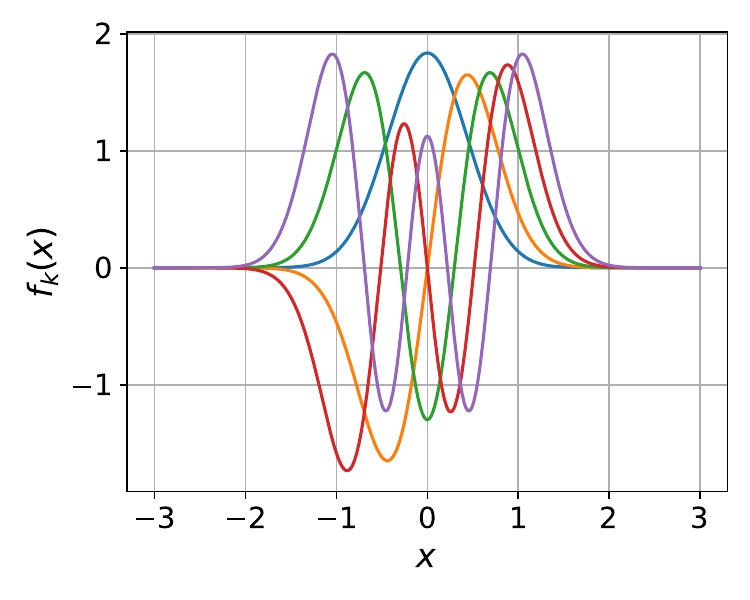}  
\includegraphics[width=0.32\linewidth]{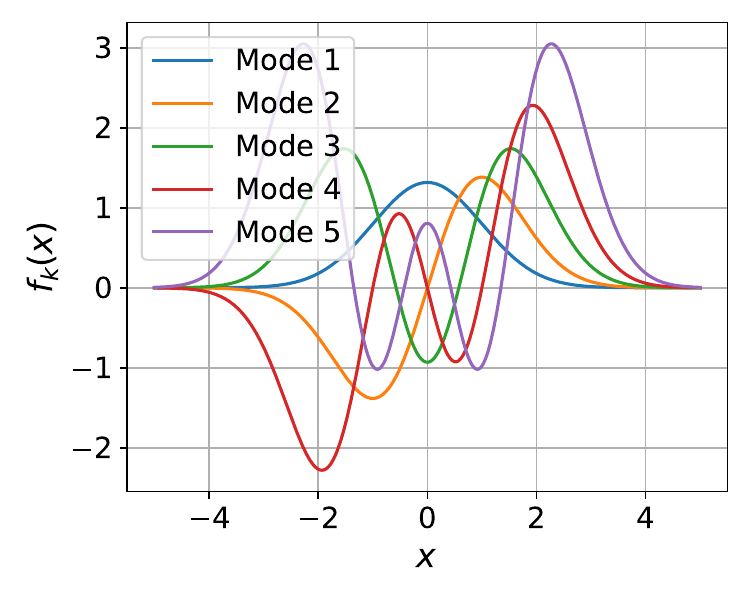}  
\includegraphics[width=0.32\linewidth]{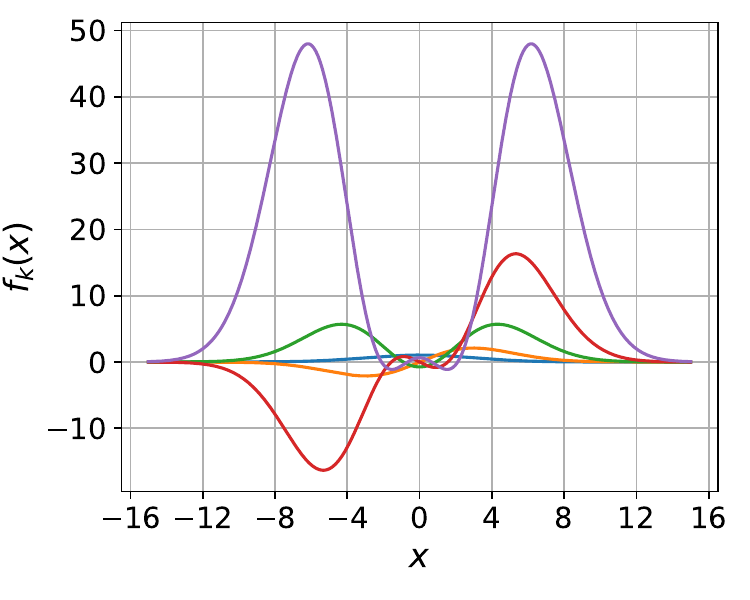}  
\end{center}
\caption{\label{fig:sqexp1D_evecs} Select eigenfunctions of a 1D KLE with a squared-exponential covariance kernel and a setup similar to Fig.~\ref{fig:sqexp1D_evals},  $\lc=\{\sigma_x/4, \sigma_x, 4\sigma_x\}$ from left to right. Note the differences between these plots for the horizontal and the vertical scales.}
\end{figure}

Figure~\ref{fig:sqexp1D_evecs} displays the first five eigenfunctions for the three cases shown in the previous figure. Note the horizontal scale differences between these cases. This scale is controlled by the $1/\sqrt{c-a}$ factor in the exponential term in Eq.~\eqref{eq:KLgauss_anl}. This factor is $\{0.6,1.4,5\}$, respectively, for the three cases presented in this figure, and the corresponding eigenfunctions consistently extend approximately $3$ times this factor on either side of $x=0$.

\subsection{Numerical Approximations via the Fredholm Equation}
\label{sec:klenum}

In practical settings we will consider a discrete mesh $\bm{X}=\{\bx_j, j=1,\ldots,\Nx\}$, where samples of the stochastic process are available or will be generated. Associated with this grid are quadrature weights $w_j$ that will be used to solve Eq.~\eqref{eq:fred} numerically. In this report we will use trapezoidal quadrature for cases where samples points are chosen at random and Gauss quadrature when one can choose grid points according to specific quadrature rules.

The covariance matrix corresponding to the covariance function evaluated on this grid will be denoted by $K$ with $K_{ij}=C(\bx_i,\bx_j)$. For brevity we will use the notation for $\fs_k$ to denote the evaluation of the eigenfunctions $f_k(\bx)$ over $\bm{X}$, with $f_k(\bm{x}_j)\rightarrow \fs_{k,j}$. 
The orthonormality of eigenfunctions, via Eq.~\eqref{eq:kleorth}, implies that numerically 
\begin{equation}
\int_{\mathfrak{D}} f_k(\bm{x}) f_l(\bm{x})\px(\bx)\dd \bx \approx \sum_{j=1}^\Nx w_j \fs_{k,j} \fs_{l,j}=\delta_{kl} .
\label{eq:evorth}
\end{equation}
The Fredholm equation (Eq.~\eqref{eq:fred}) is discretized as follows
\begin{equation}
\int_{\mathfrak{D}}C(\bx_1,\bx_2) f(\bx_2)\px(\bx_2)\dd \bx_2
\approx \sum_{j=1}^\Nx w_j \Ks_{ij} \fs_j = \lambda \fs_i,\quad i=1,\ldots,\Nx.
\label{eq:fredholm}
\end{equation}
In matrix notation
\begin{equation}
\Ks\,\Ws\, \fs = \lambda \fs
\label{eq:es1}
\end{equation}
where $\Ws$ is a diagonal matrix with the quadrature weights on the diagonal, $\Ws_{jj}=w_j$. In this form, the $\lambda$ and $\fs$ are eigenvalues and eigenvectors of $\Ks\,\Ws$. However, unless $\Ws\propto I$, it is clear that $\Ks\,\Ws$ is not symmetric. As a result, we do not have the desired symmetric positive definite eigenproblem with the above decomposition in terms of real eigenvalues and eigenvectors.

This can be alleviated by rewriting Eq.~\eqref{eq:es1} as follows
\begin{eqnarray}
\Ks\Ws\, \fs&=&\lambda \fs \notag\\
(\Ks\Ws^{1/2})\,(\Ws^{1/2}f)&=& \lambda \fs \notag\\
\underbrace{(\Ws^{1/2}\Ks\Ws^{1/2})}_{\As}\, \underbrace{(\Ws^{1/2}\fs)}_{\hs}&=& 
\lambda\, \underbrace{\Ws^{1/2}\fs}_{\hs}\notag\\
\As \hs &=&\lambda \hs \label{eq:esmod}
\end{eqnarray} 
The eigenvectors of $\As:=\Ws^{1/2}\Ks\Ws^{1/2}$ are orthonormal
\begin{equation}
\sum_{j=1}^m \hs_{k,j}\hs_{l,j}=\delta_{kl}, \label{eq:orthonorm}    
\end{equation}
and, since $\hs=\Ws^{1/2}\fs$, this results in eigenvectors $\fs$ satisfying Eq.~\eqref{eq:evorth}.

Once eigenvalues $\lambda_k$ and eigenvectors $\fs_k$ are available, one can obtain samples of uncorrelated RVs using a discretized form of Eq~\eqref{eq:rvest} 
\begin{equation}
\xi_k^{(p)}=\frac{1}{\sqrt{\lambda_k}}\sum_{j=1}^m w_j Y_j^{(p)}\fs_{k,j}\,\,\textrm{and}\,\,
\zeta_k^{(p)}=\sum_{j=1}^m w_j Y_j^{(p)}\fs_{k,j}
\label{eq:xis1}
\end{equation}
where $Y_j^{(p)}$ is the $p-$th sample at location $\bm{x}_j$ of stochastic process $Y$, and $\xi_k^{(p)}$ is the $p-$th sample of RV $\xi_k$.

\subsection{Numerical Approximation via Singular Value Decomposition}
\label{sec:svd}

Let $S$ be the matrix constructed with samples of $Y$. Each column $i$ of $S$ is one sample $Y_i$. Each row $j$ corresponds to grid point $\bx_j\in\mathfrak{D}$. Let $m$ be the number of grid points and $n$ the number of samples. The covariance matrix $K$ can be constructed from samples $S$ as
\begin{equation}
\stackunder{\Ks}{\sss m\times m} =  \frac{1}{n-1}\stackunder{\Ss}{\sss m\times n}\,  
\stackunder{\Ss^T}{\sss n\times m}
\label{eq:covspl}
\end{equation}
Please note that in the above equation and those following below we display the size of matrices to highlight the equivalence between the Singular Value Decomposition (SVD) algorithm and the discretized Fredholm equation approach described in the previous section. The SVD algorithm avoids the estimation of the covariance matrix, as follows.

Consider the sample matrix $\Ss$ scaled as in $\Ss/\sqrt{n-1}$. Applying SVD to this matrix gives
\begin{equation}
\frac{1}{\sqrt{n-1}}\,\stackunder{\Ss}{\sss m\times n} =  \stackunder{\Us}{\sss m\times m} 
\, \stackunder{\SGs}{\sss m\times n}\,\stackunder{\Vs^T}{\sss n\times n} 
\label{eq:svddec}
\end{equation} 
where $\Us$ and $\Vs$ are unitary matrices and $\SGs$ is a rectangular matrix, generally non-square, with non-negative singular values $\sigma_i$ on the diagonal. Let's follow with some algebraic manipulations.
\begin{eqnarray}
\frac{1}{\sqrt{n-1}} \Ss\, \frac{1}{\sqrt{n-1}} \Ss^T=\Ks&=&\Us\,\SGs\, \underbrace{\Vs^T\,
   \Vs}_{\mathbb{I}}\,\SGs^T\, \Us^T\notag \\
\Ks&=&\Us\,\underbrace{\SGs\,\SGs^T}_{\SGs^2}\, U^T, \qquad\text{where}\quad \SGs^2 \in \mathbb{R}^{m\times m}  \notag \\
\Ks&=&\Us\,\SGs^2 \, \Us^T  \\
\Ks\, \Us&=&\Us\,\SGs^2
\end{eqnarray}
The last expression above is the eigensystem for the covariance matrix $\Ks$ when the quadrature weights are all equal to one, i.e. $W=\mathbb{I}$. For this case the eigenvalues $\lambda_i$ of the covariance matrix are the squares of singular values $\sigma_i$ of the scaled sample matrix, $\Ss/\sqrt{n-1}$. The left singular vectors of $\Ss/\sqrt{n-1}$, that form the columns of $\Us$, are the right eigenvectors of $\Ks$.

We can adapt the scaling of $\Ss$ to accommodate non-uniform quadrature weights. Starting with the eigensystem in Eq.~\eqref{eq:esmod}, $\As\hs=\lambda \hs$, we rewrite $\As$ as
\begin{equation}
\As=\Ws^{1/2}\, \Ks\, \Ws^{1/2}=\Ws^{1/2}\,\frac{1}{n-1}\Ss\Ss^T\, \Ws^{1/2}
=\left(\frac{1}{\sqrt{n-1}}\Ws^{1/2}\Ss\right)\, \left(\frac{1}{\sqrt{n-1}}\Ws^{1/2}\Ss\right)^T
\end{equation}
Therefore the eigenvalues of $A$ (corresponding to the discretization of the Fredholm equation) are the squares of the singular values of the weighted-scaled sample matrix:
\begin{equation}
\frac{1}{\sqrt{n-1}}\Ws^{1/2}\Ss
\label{eq:scaleW}
\end{equation}

\subsection{Uniform Scaling of Quadrature Weights}
\label{sec:uniw}

Consider a setup with where the covariance matrix $\Ks$ is scaled by a constant factor, e.g. $\Ks\rightarrow \sigma^2 \Ks$, i.e. the stochastic process has variance $\sigma^2$. The eigensystem in Eq.~\eqref{eq:esmod} becomes
\begin{equation}
\underbrace{(\sigma\Ws^{1/2})\, \Ks\, (\sigma\Ws^{1/2})}_{\As'}\, \underbrace{\Ws^{1/2}\, \fs'}_{\hs'}=
\lambda'\, \underbrace{\Ws^{1/2}\,\fs'}_{\hs'}.
\label{eq:kleug}
\end{equation}
The scaling of $\Ks$, e.g. increasing or decreasing the variance of the random process, is equivalent to a uniform scaling of the symmetric matrix $\As$. Given that the eigenvectors are numerically estimated using algorithms that yield orthonormal vectors, this is equivalent to a new set of eigenvalues that are similarly scaled compared to the original system, $\lambda'=\sigma^2\lambda$, and the same set of orthonormal eigenvectors, $\hs'=\hs$. Thus
\[
\fs'_k = \fs_k = \Ws^{-1/2}\hs_k.
\]

\section{Numerical Examples}
\label{sec:numerical}
We provide one-dimensional (1D), as well as 2D and 3D examples that illustrate the numerical aspects associated with Karhunen-Lo\`{e}ve expansions.

\subsection{Uniform Spatial Density with Exponential Covariance Kernel}
Consider a 1D configuration, with the spatial coordinate $x$ corresponding to a uniform density on $\mathfrak{D}=[0,1]$, resulting in $\px(x)=1$. The spatial domain is discretized by $\Nx$ uniformly spaced cell-centered grids
\begin{equation}
  x_i = \frac{i - \tfrac{1}{2}}{\Nx}, 
  \qquad w_i = \Delta x = \frac{1}{\Nx}, \quad i = 1,\ldots,\Nx,
\end{equation}
giving equal quadrature weights. Three grid resolutions are considered, $\Nx \in \{64, 256, 1024\}$, in the examples presented in this section. We consider a zero-mean, unit-variance Gaussian random field $Y(x,\omega)$ on the interval $\mathfrak{D}$ with an exponential covariance kernel given by Eq.~\eqref{eq:expcov}, parametrized by the correlation length $\lc > 0$. 

\begin{figure}[htb!]
\centering
\begin{subfigure}[b]{0.32\textwidth}
\includegraphics[width=\textwidth]{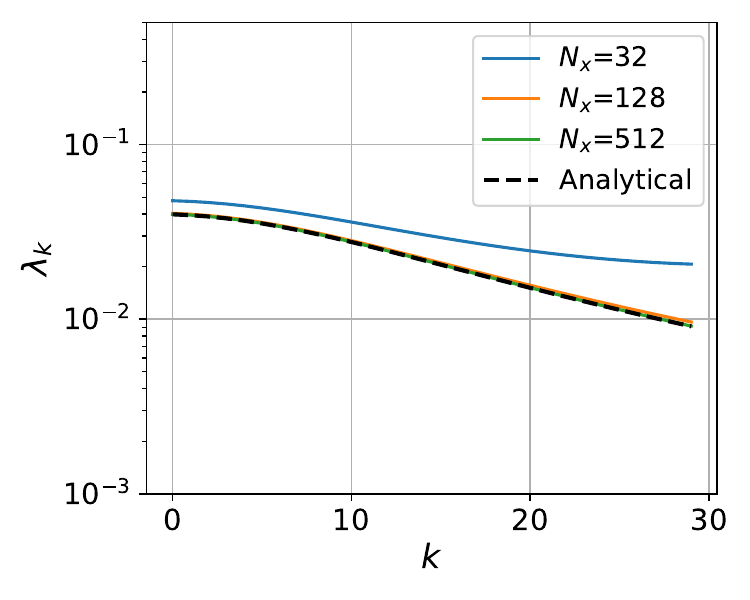}
\caption{$\lc = 0.02$}
\label{fig:evals_lc0}
\end{subfigure}
\hfill
\begin{subfigure}[b]{0.32\textwidth}
\includegraphics[width=\textwidth]{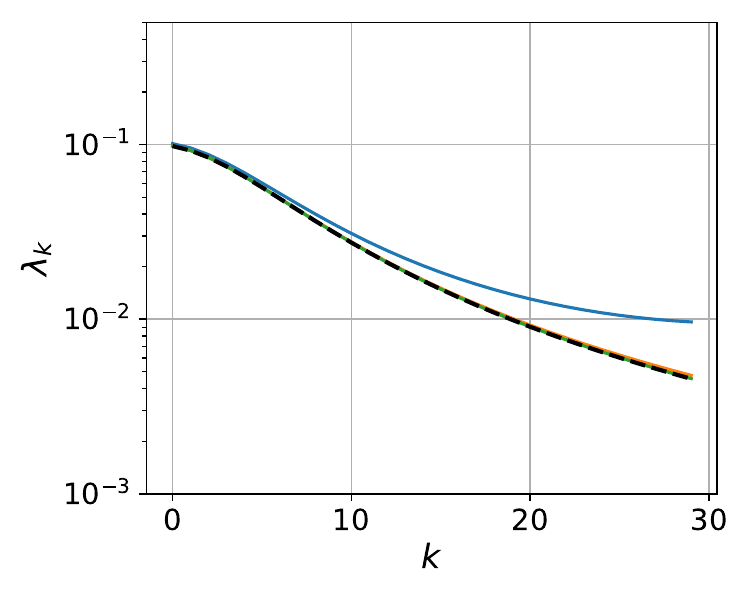}
\caption{$\lc = 0.05$}
\label{fig:evals_lc1}
\end{subfigure}
\hfill
\begin{subfigure}[b]{0.32\textwidth}
\includegraphics[width=\textwidth]{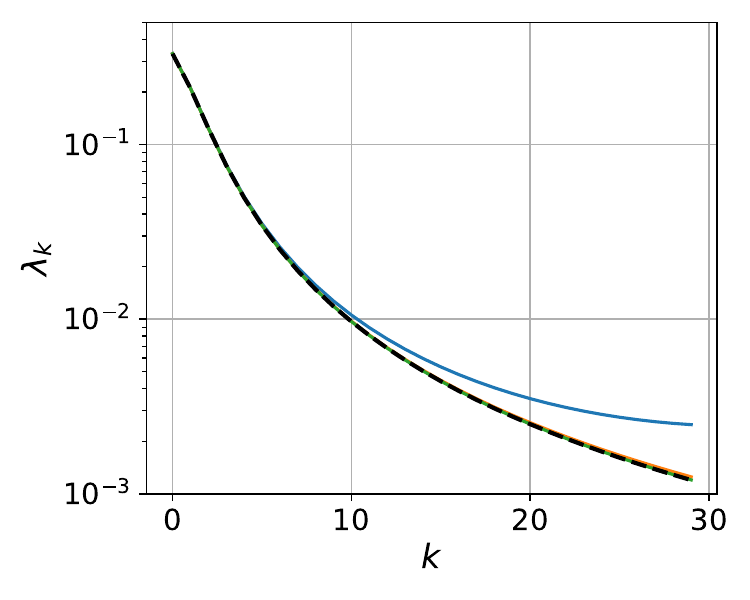}
\caption{$\lc = 0.2$}
\label{fig:evals_lc2}
\end{subfigure}
\caption{Numerical eigenvalues (solid lines, $\Nx \in \{32,128,512\}$) versus analytical solution (dashed line) for each correlation length.}
\label{fig:ev_unif_fred}
\end{figure}

Results for three correlation lengths are examined below, $\lc \in \{0.02, 0.05, 0.2\}$. Figure~\ref{fig:ev_unif_fred} compares the first 30 numerical eigenvalues against the analytical solution, obtained through expressions presented in Sec.~\ref{sec:klexpanl}, for each correlation length. The numerical results presented in this figure were obtained by first computing the covariance matrix corresponding to each choice of computational grids, followed by the numerical solution of the Fredholm equation~\eqref{eq:fred}.
The eigenvalues decay more slowly for shorter $\lc$, reflecting the fact that more terms are required in the KL expansion to represent rougher realizations produced by the exponential kernel over small scales. For the coarse grid, $\Nx=32$, there are fewer than one grid point per correlation length for the case with $\lc=0.02$ leading to visible discrepancies between the numerical and analytical results. As the correlation length increases the coarse grid solution approaches the analytical results at lower frequencies. For settings with larger grids, $\Nx=128$ and $\Nx=512$, the grid density is sufficient even for the shortest correlation length leading to a good agreement between numerical and analytical results.

\begin{figure}[h]
  \centering
  \includegraphics[width=\textwidth]{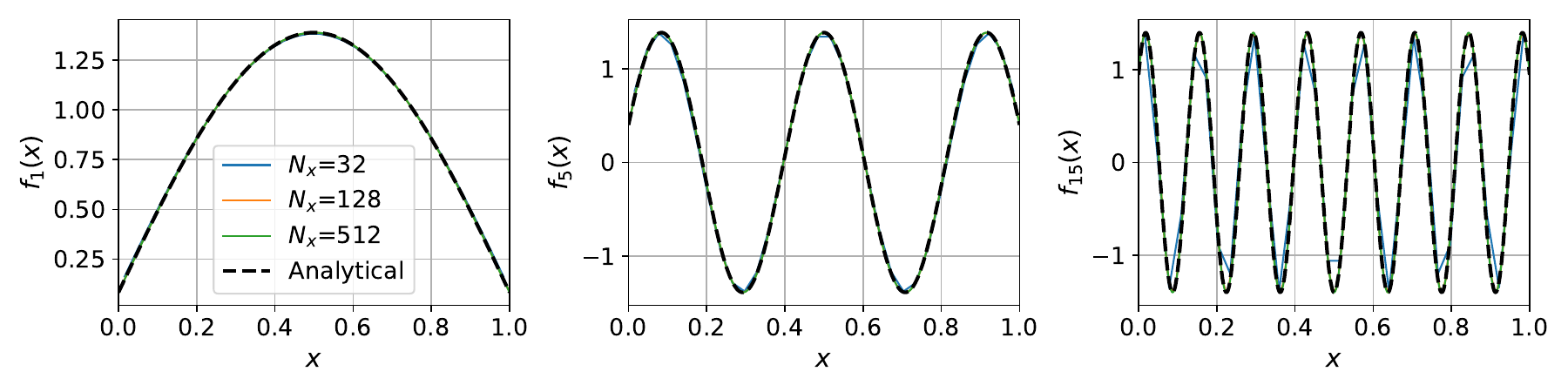}
  \caption{Eigenfunctions $1$, $5$, and $15$ for $\lc = 0.02$.}
  \label{fig:emodes_unif_lc0}
\end{figure}

Figures~\ref{fig:emodes_unif_lc0} shows eigenfunctions 1, 5, and 15, for $\lc=0.02$.  The numerical solutions (solid color lines) were sign-aligned with the corresponding analytical eigenfunctions (dashed lines). Unlike the eigenvalue spectrum, the eigenfunctions are less sensitive to discretization errors and the results are close to analytical results even for the coarse grid case, $\Nx=32$. For higher frequencies, the grid effects start to be noticed as coarse grids do not capture the oscillatory structure at certain locations.

\begin{figure}[htb!]
\centering
\includegraphics[width=0.9\textwidth]{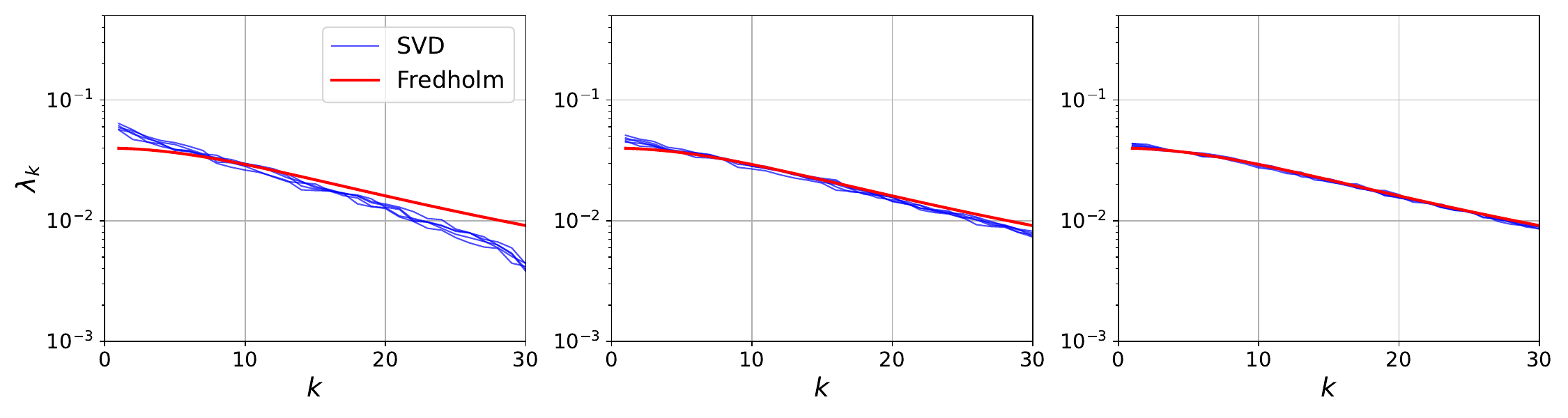}
\caption{SVD eigenvalues (blue lines, 5 seeds) vs.\ Fredholm reference (red lines) for $\lc = 0.02$. The SVD results in the left, middle, and right frames correspond to $\Ns=128$, $512$, and $2048$, respectively.}
\label{fig:eval_svd_lc0}
\end{figure}

\begin{figure}[htb!]
\centering
\includegraphics[width=0.9\textwidth]{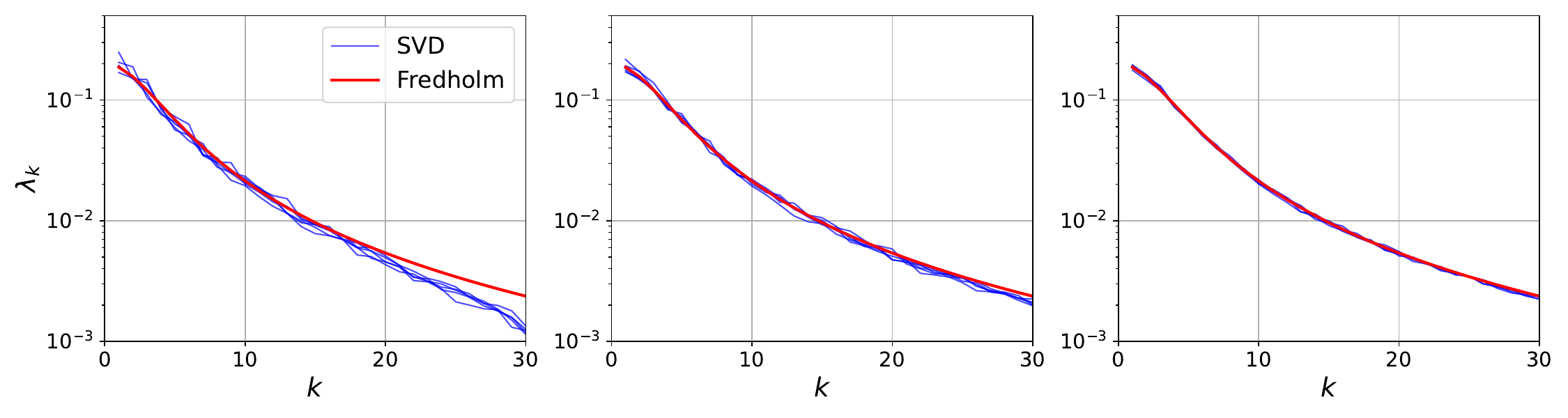}
\caption{SVD eigenvalues (blue lines, 5 seeds) vs.\ Fredholm reference (red lines) for $\lc = 0.1$. SVD results are based on the same $\Ns$ as in Figure~\ref{fig:eval_svd_lc0}.}
\label{fig:eval_svd_lc2}
\end{figure}

In practice, the covariance kernel is typically not known analytically and must be
estimated from an ensemble of random field realizations. Below we examine the effect of the number of samples on the eigenvalue spectra and compare the SVD and Fredholm approaches to estimate eigenvalues and eigenfunctions for KL expansions. For this set of results we employ a fixed, uniform, grid with $\Nx = 1024$ on $\mathfrak{D}=[0,1]$ and use the following steps to assemble the numerical tests.
\begin{enumerate}
\item for each choice of correlation length $\lc$ we assemble the covariance matrix via Eq.~\eqref{eq:expcov}, then use the Fredholm KLE to compute the first $\Nm$ eigenpairs $(\lambda_k,\fs_k)$. Given the dense grid employed here, the numerical solution of the Fredholm KLE is close to the analytical solution, per the comparison presented above in Figs.~\ref{fig:ev_unif_fred} and~\ref{fig:emodes_unif_lc0}.
\item draw $\Ns$ standard-normal coefficients for each mode $k$, $\xi^{(j)}_k,\, k=1,\ldots,\Nm$ and form samples
\begin{equation}
  Y^{(j)} = \sum_{k=1}^\Nm \xi_k^{(j)}\sqrt{\lambda_k}\,\fs_k.
\end{equation}
\item Run the SVD algorithm presented in Sec~\ref{sec:svd} on the centered sample matrix $S$ assembled with the samples $Y$ to recover $\tilde\lambda_k$ and the empirical eigenvectors.
\end{enumerate}
These steps are repeated for 5 independent random seeds for each choice of $\lc \in \{0.02, 0.1\}$ and $\Ns \in \{128, 512, 2048\}$. The SVD eigenvalues (blue lines, one curve per seed) are compared against the Fredholm reference (red line) in Figs.~\ref{fig:eval_svd_lc0} and~\ref{fig:eval_svd_lc2}. Variability across seeds decreases as the number of samples grows, bringing the SVD spectra close to the Fredholm eigenvalues across all modes shown here. While not immediately evident from these semi-log plots, short correlation lengths ($\lc = 0.02$) require a larger $\Ns$ to resolve the slower-decaying tail.

\begin{figure}[htb!]
\centering
\includegraphics[width=0.9\textwidth]{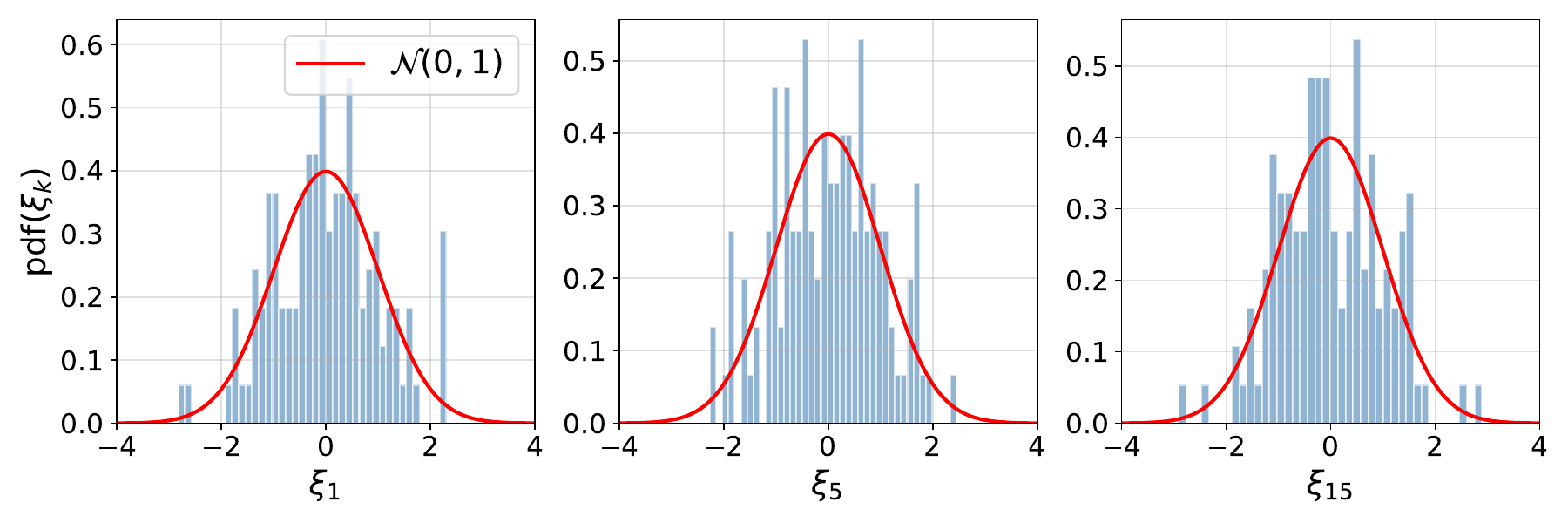}
\caption{Histograms of select KLE coefficients for $\lc=0.1$ and $\Ns=128$.}
\label{fig:svd_hist_n128}
\end{figure}

\begin{figure}[htb!]
\centering
\includegraphics[width=0.9\textwidth]{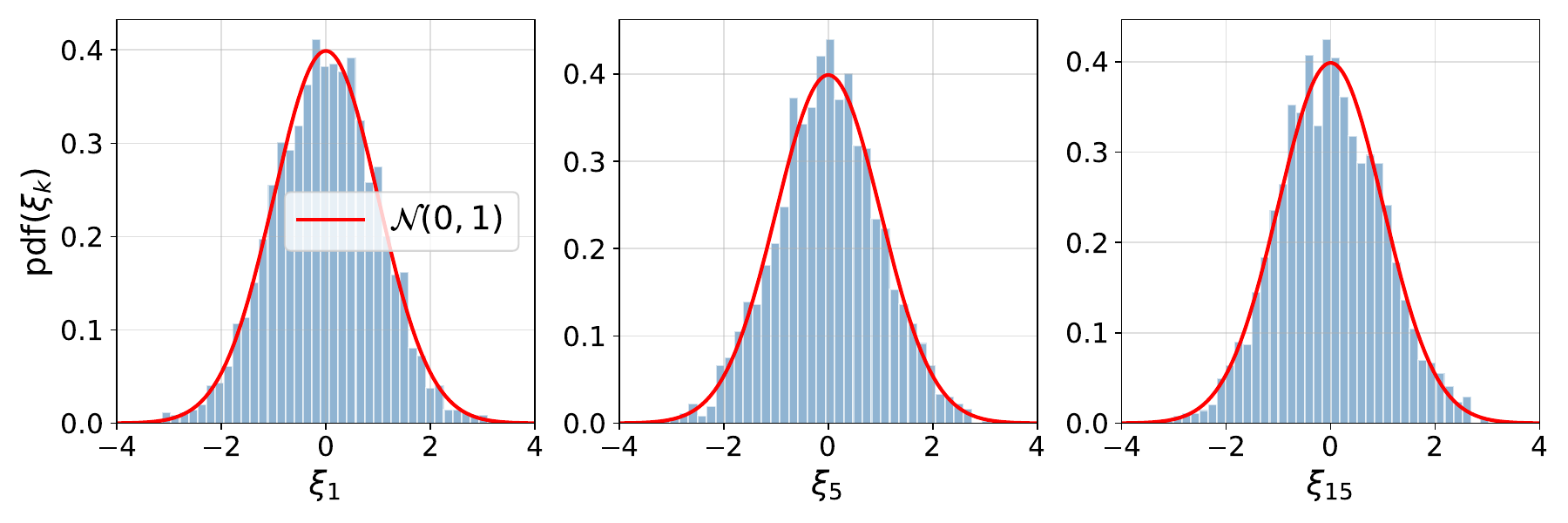}
\caption{Histograms of select KLE coefficients for $\lc=0.1$ and $\Ns=2048$.}
\label{fig:svd_hist_n2048}
\end{figure}

The KL coefficients $\xi_k$ extracted from the SVD expansion should be approximately $\mathcal{N}(0,1)$ given that the underlying generative process is Gaussian. Histograms of $\xi_1$, $\xi_5$, and $\xi_{15}$ are shown in Figs.~\ref{fig:svd_hist_n128} and~\ref{fig:svd_hist_n2048} for $\lc = 0.1$ and two sample counts.
The agreement with $\mathcal{N}(0,1)$ improves with more samples, as expected.

To quantify convergence of the marginal distributions to $\mathcal{N}(0,1)$ we compute the Kullback-Leibler (KL) divergence
\begin{equation}
  D_{\mathrm{KL}}\!\left(p\,\big\|\,\hat{q}_k^{(\Ns)}\right)
  = \int_{-\infty}^{\infty}\ln\left(\frac{p(\xi)}{\hat{q}_k^{(\Ns)}(\xi)}\right)\,p(\xi)\dd\xi,
\end{equation}
where $p = \mathcal{N}(0,1)$ is the theoretical standard normal solution and $\hat{q}_k^{(\Ns)}$ is the kernel-density estimate (KDE) of $p(\xi_k)$ built from $\Ns$ samples. These results use a correlation length $\lc = 0.1$, with \hbox{$\Ns \in \{32, 128, 512, 2048\}$}, and $100$ independent random seeds.

\begin{figure}[htb!]
\centering
\includegraphics[width=0.9\textwidth]{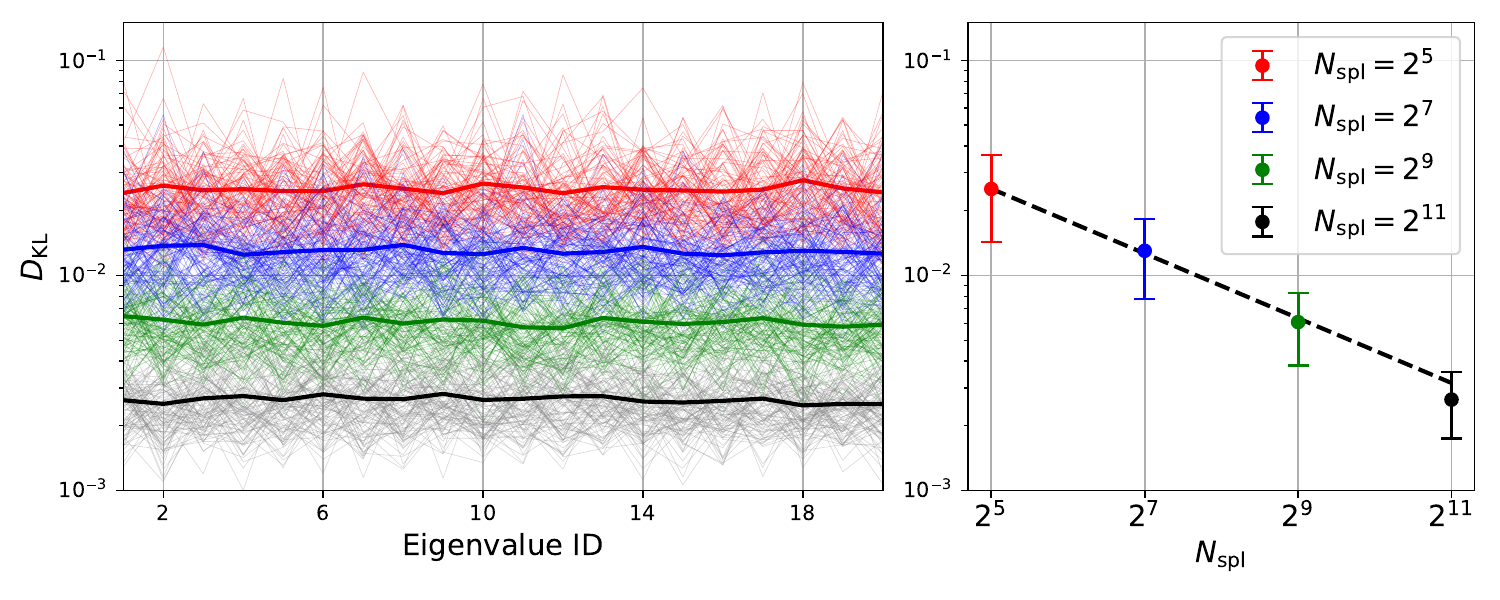}
\caption{KL divergence of empirical KLE coefficient densities from $\mathcal{N}(0,1)$ for $\lc = 0.1$.
\emph{Left:} $D_{\mathrm{KL}}$ vs.\ mode index; thin lines are individual samples, bold lines are sample means.
\emph{Right:} mean $\pm$ std over all modes and seeds vs.\ $\Ns$; dashed line is a the fitted $1/\sqrt{n}$ trendline.}
\label{fig:kldiv}
\end{figure}

The left panel of Fig.~\ref{fig:kldiv} shows $D_{\mathrm{KL}}$ as a function of mode index for each sample count; thin lines represent $D_{\mathrm{KL}}$ for individual seeds and bold lines are corresponding sample means. The right panel shows mean $\pm$ one standard deviation over all modes and seeds as a function of $\Ns$, together with a fitted $1/\sqrt{n}$ trendline via least-squares. The KL divergence decreases monotonically with $\Ns$ across all modes, confirming that the SVD algorithm recovers the correct statistics for the KL coefficients as the sample size grows. The $1/\sqrt{n}$ trend is consistent with the expected Monte Carlo convergence for the sampling-based approach.

\subsection{Normal Spatial Distribution with Squared-Exponential Covariance}

This section explores the numerical solution of Eq.~\eqref{eq:fredholm} for the one-dimensional squared-exponential covariance kernel, in Eq.~\eqref{eq:sqexpcov}, on the infinite real line, with the Gaussian distribution $\px(x) = \mathcal{N}(0,\sigma_x^2)$. The numerical results presented here employ two classes of quadrature grids
\begin{enumerate}
\item \textbf{Monte Carlo:} grid nodes drawn independently from $\mathcal{N}(0,\sigma_x^2)$, with weights proportional to the Gaussian PDF evaluated at each node.
\item \textbf{Gauss-Hermite (GH):} quadrature nodes and weights optimal for the Gaussian distribution, constructed via Gauss-Hermite quadrature, with a change of variable to adapt the quadrature points and weights from $\px(x)\propto \exp(-x^2)$ to $\px(x)\propto \exp(-x^2/(2\sigma_x^2))/\sigma_x$. 
\end{enumerate}
Three correlation lengths $\lc \in \{\sigma_x/4,\,\sigma_x,\,4\sigma_x\}$ are examined. Results presented in this section remain qualitatively similar as they depend up to floating point precision on the ratio $\rho=\lc/\sigma_x$.

\begin{figure}[h]
  \centering
  \includegraphics[width=0.85\linewidth]{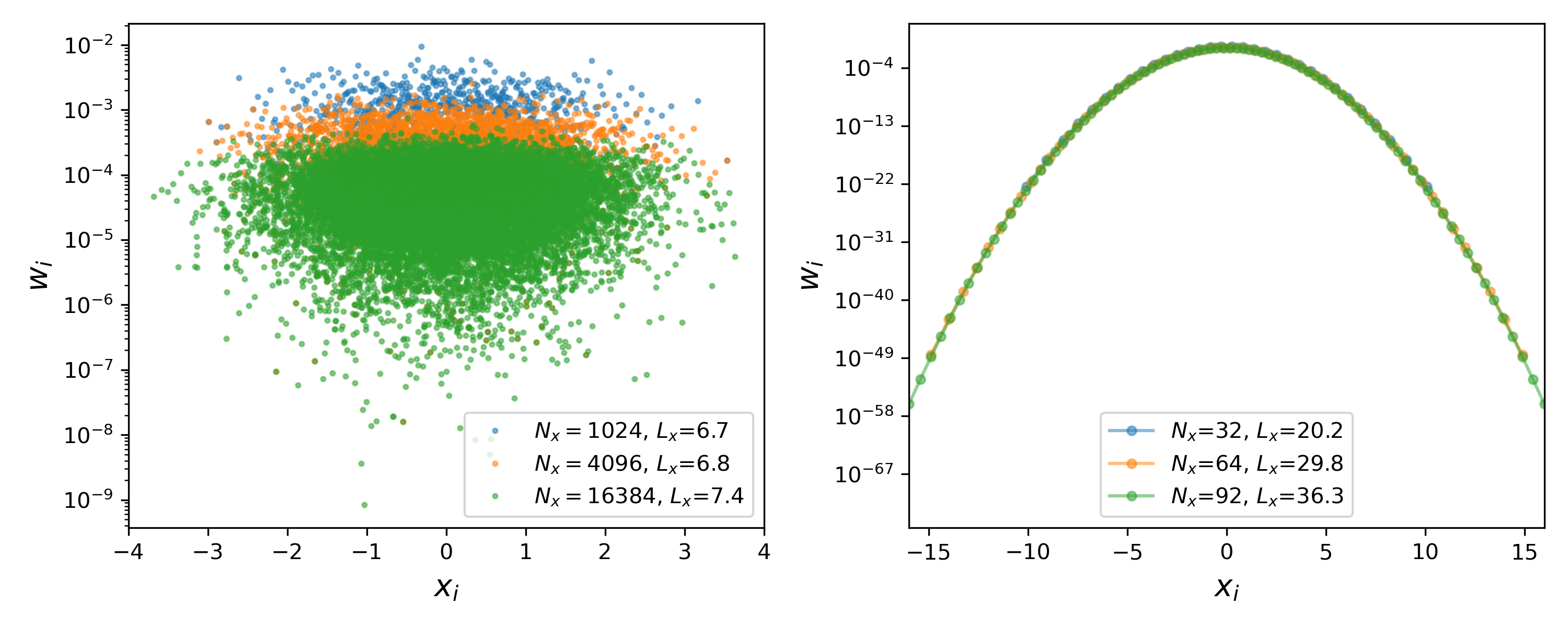}
  \caption{Quadrature nodes $x_i$ and weights $w_i$ for MC Gaussian grids (left) and Gauss-Hermite grids (right), for $\sigma_x=1$ and $\lc=\sigma_x/4$.}
  \label{fig:grids}
\end{figure}

Figure~\ref{fig:grids} compares the two grid types in terms of nodes' location $x_i$ and associated quadrature weights $w_i$, for $\sigma_x=1$ and $\ell_c=\sigma_x/4$.
The MC grids, in the left panel, place points within roughly $[-4\sigma_x, 4\sigma_x]$, with weights that scatter over approximately 5 orders orders of magnitude.  This large weight spread reflects the non-uniform density of a random Gaussian draw.  The effective domain extent grows only negligibly with $\Nx$. The domain lengths $L_x$ are measured as the distance between the first and last grid points and are approximately $6.7$, $6.8$, and $7.4$ for $\Nx \in\{2^{10}$, $2^{12}$, $2^{14}\}$.

The GH grids, in the right panel, span a much wider range, reaching $\pm 15\sigma_x$ for $\Nx = 92$, indicating a wider coverage of the tails of the distribution. Their weights follow the Gaussian PDF exactly by construction and cover nearly 55 orders of magnitude, giving machine-precision integration of all polynomials up to degree $2\Nx-1$.

Below, we will compare numerical results obtained via quadrature, either MC or GH, with analytical results presented in Sec.~\ref{sec:anlsqexp}. The covariance matrix is evaluated via Eq.~\eqref{eq:sqexpcov} for each set of grids and the KLE eigenvalues and eigenvectors are computed as the numerical solution of the discretized Fredholm equation~\eqref{eq:fredholm}.

\begin{figure}[!htb]
\begin{center}
\includegraphics[width=0.85\linewidth]{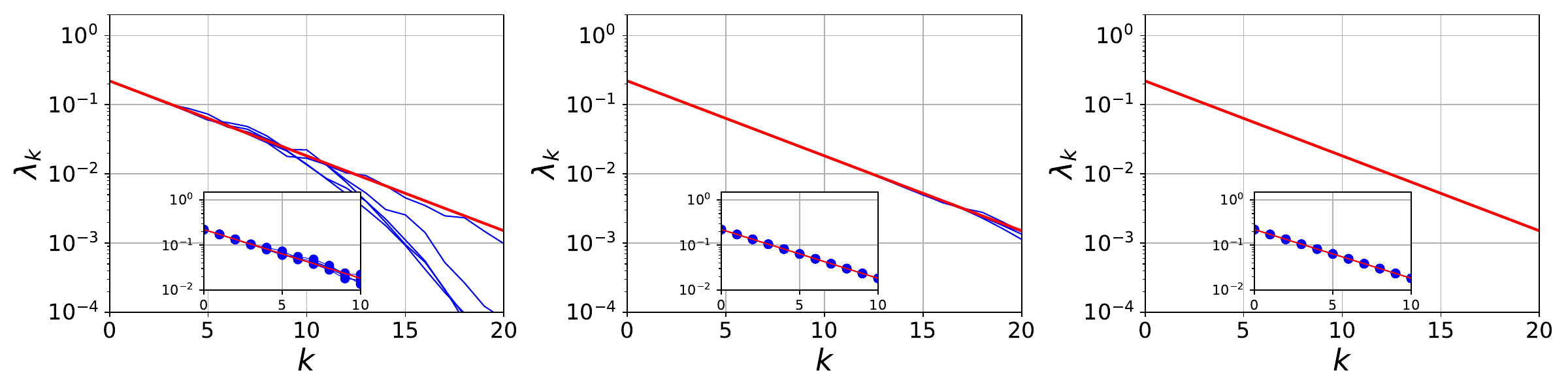}  
\includegraphics[width=0.85\linewidth]{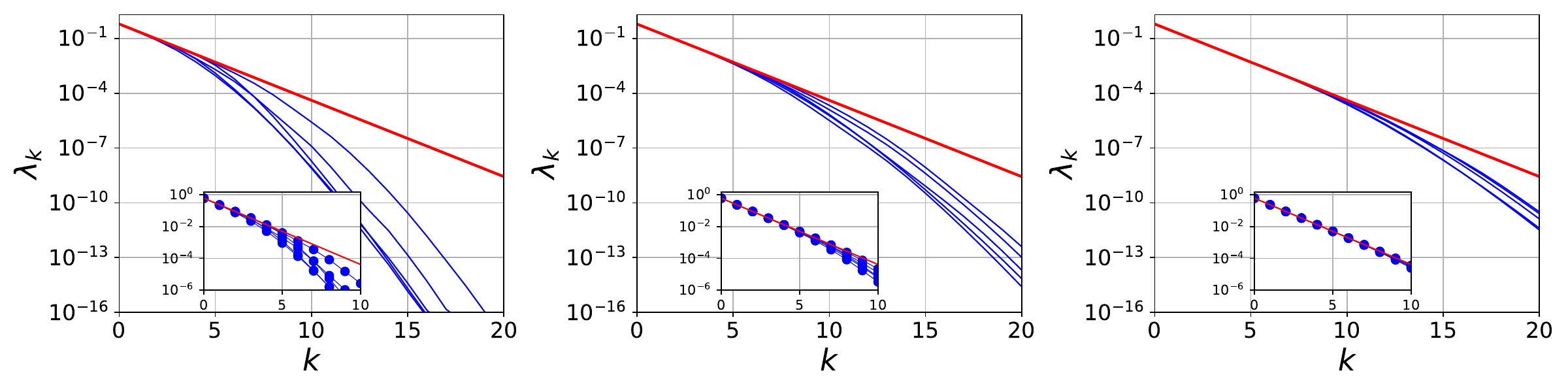}  
\includegraphics[width=0.85\linewidth]{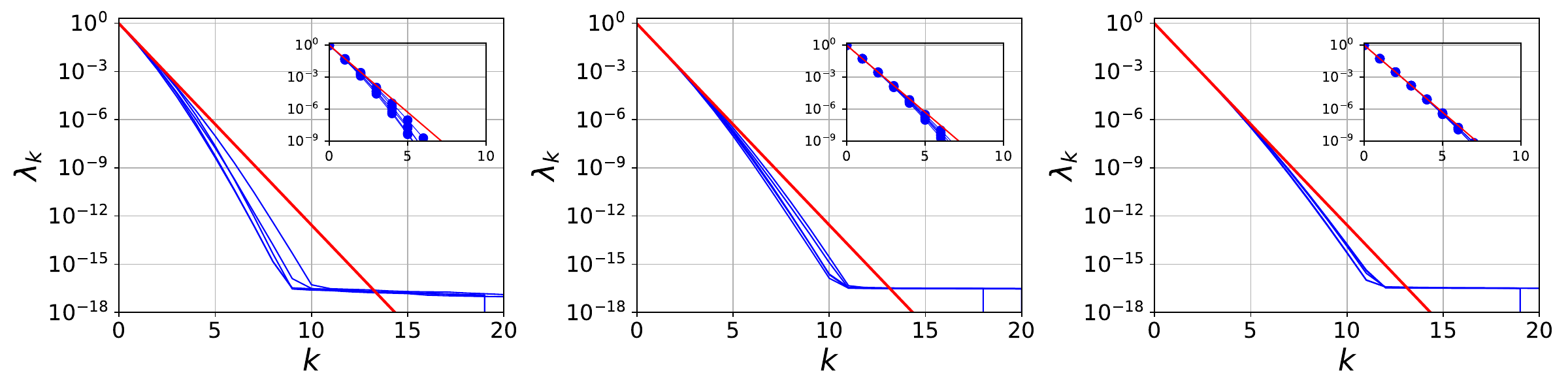}  
\end{center}
\caption{\label{fig:ev_gauss_mc} Karhunen-Lo\`{e}ve expansion eigenvalues for correlation lengths $\lc=\{\sigma_x/4,\sigma_x,4\sigma_x\}$ from top to bottom rows, respectively. Frames correspond to grid sizes $\Nx\in\{2^6,2^{10},2^{14}\}$ from left to right. The blue lines show eigenvalues for different MC grid realizations, while the red lines correspond to analytical values obtained via Eq.~\eqref{eq:eval_anl_expression}.}
\end{figure}

Fig.~\ref{fig:ev_gauss_mc} shows eigenvalue spectra for cases with different correlation lengths $\ell_c$ and number of grid points $\Nx$. 
An inset zooms into the first ten modes. For the shorter correlation length (top row), the leading 5-6 eigenvalues are captured well even for a small grid size ($\Nx=2^6$). Increasing the number of grid points to $2^{10}$ and beyond, leads to numerical eigenvalue spectra close to analytical results for the range of modes shown in this figure. For $\lc = \sigma$ the spectrum drops quadratically (middle row). At $\Nx = 2^6$, the numerical eigenvalues begin to underestimate the analytical values after $k \approx 4$; by $k \approx 10$ they scatter widely over several orders of magnitude, and the individual seed curves diverge completely. For larger grid sizes convergence improves and the numerical and analytical results agree well to twice and three times the number of modes, respectively. For $\lc = 4\sigma_x$ the eigenvalue decay is such that by  $k \approx 15$ the eigenvalues hit the double-precision noise floor ($\approx 10^{-16}$) regardless of the grid size and MC instance. Somewhat counterintuitively, all grid sizes diverge from the analytical results after 3-5 modes, with lower $\Nx$ cases diverging earlier.
Nevertheless, at this correlation length, most of the process variance, $99.9\%$, is recovered by the first three modes.

These results suggest two regimes for MC convergence.  For short $\lc$, the random Gaussian grids extend sufficiently to capture the eigenspectra of the random process.  For long $\lc$, the eigenfunctions extend well beyond the typical MC sample domain; however, because only a few eigenvalues are numerically significant, this may not be consequential. 

\begin{figure}[htb]
\begin{center}
\includegraphics[width=0.85\linewidth]{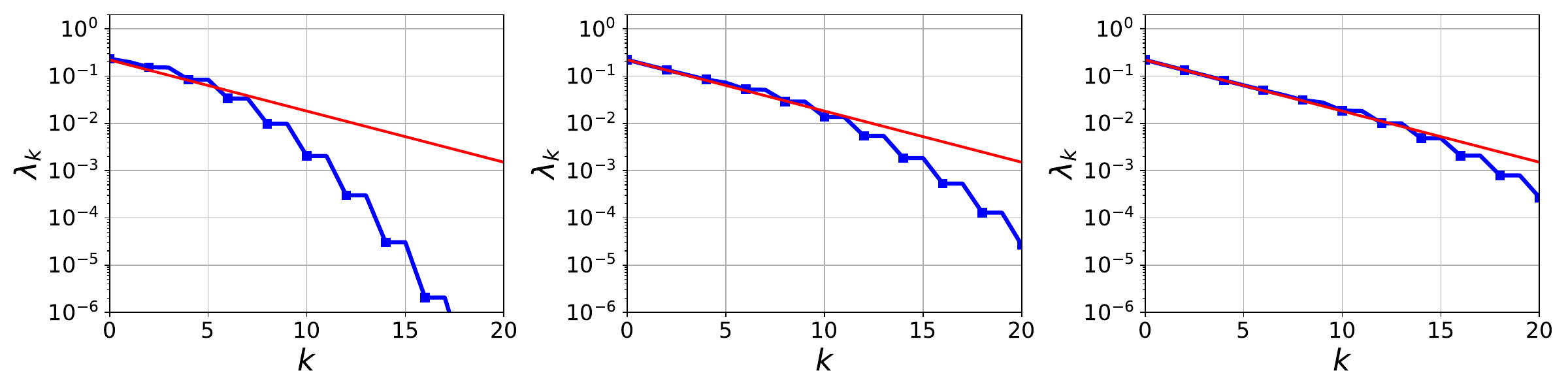}  
\includegraphics[width=0.85\linewidth]{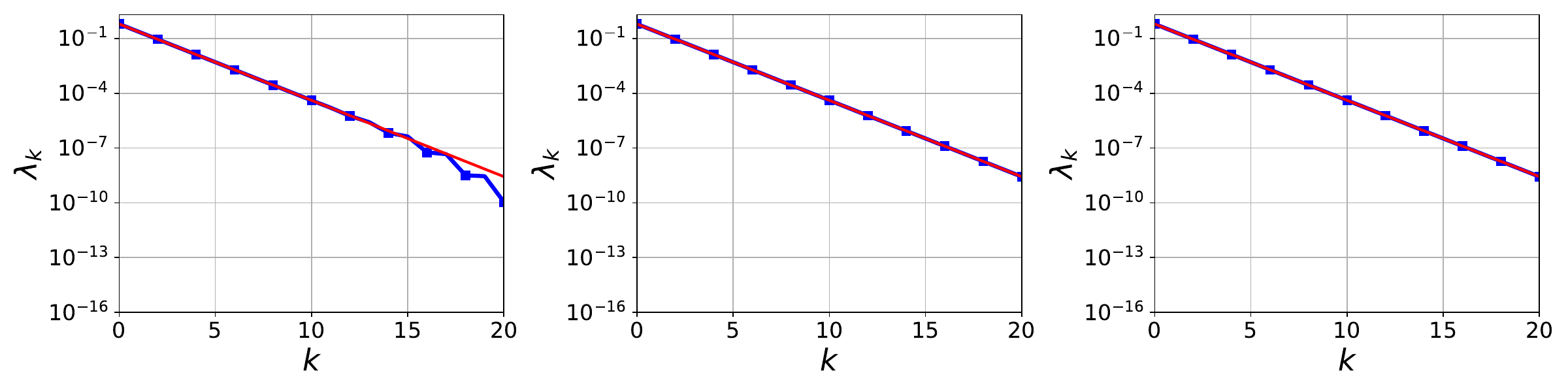}  
\includegraphics[width=0.85\linewidth]{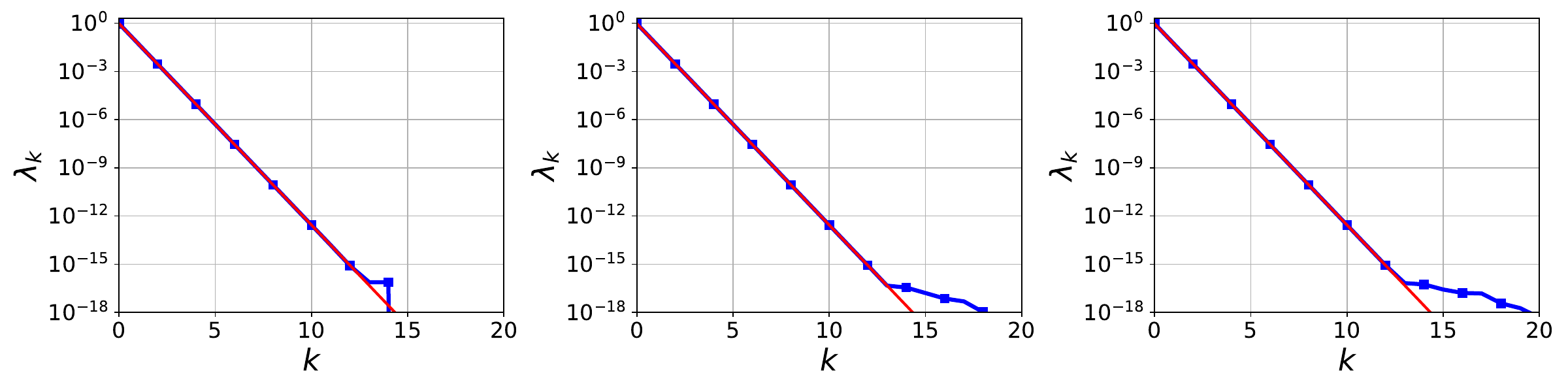}  
\end{center}
\caption{\label{fig:ev_gauss_her} Karhunen-Lo\`{e}ve expansion eigenvalues for a squared-exponential covariance matrix with the same set of correlation lengths as in Fig.~\ref{fig:ev_gauss_mc}. Frames correspond to increasing GH quadrature order from left to right $\Nx\in\{32,64,92\}$. The blue lines show eigenvalues obtained numerically via GH quadrature, while the red lines correspond to analytical values.}
\end{figure}

Fig.~\ref{fig:ev_gauss_her} shows a comparison similar to the previous figure, except now we use Gauss-Hermite (GH) quadrature. 
The error trends observed for the GH results are somewhat opposite to the MC trends. For GH applied to the short correlation length case (the first row), increased quadrature order improves the quality of the results but the agreement with the theoretical results is poorer compared to the MC cases using $\Nx=2^{10}$ or larger grids. This is a result of the MC grids concentration in regions of high values for $\px(x)$ and, for small $\lc$ values, that restricts most of the random process variance to that region. For this case, the GH results, with samples and associated weights generated based on a lengthscale $\sigma_x=4\lc$, are less accurate compared to the tightly packed MC samples. For larger correlation lengths (middle and bottom rows in Fig.~\ref{fig:ev_gauss_her}), the GH agreement with the analytical results is clearly better than the corresponding results for MC quadrature shown in the previous figure. For these cases the GH quadrature grids that are present in regions that extend in the tails of the spatial distributions, capture the random process accurately, unlike the MC approach.

\begin{figure}[!htb]
\begin{center}
\includegraphics[width=0.85\linewidth]{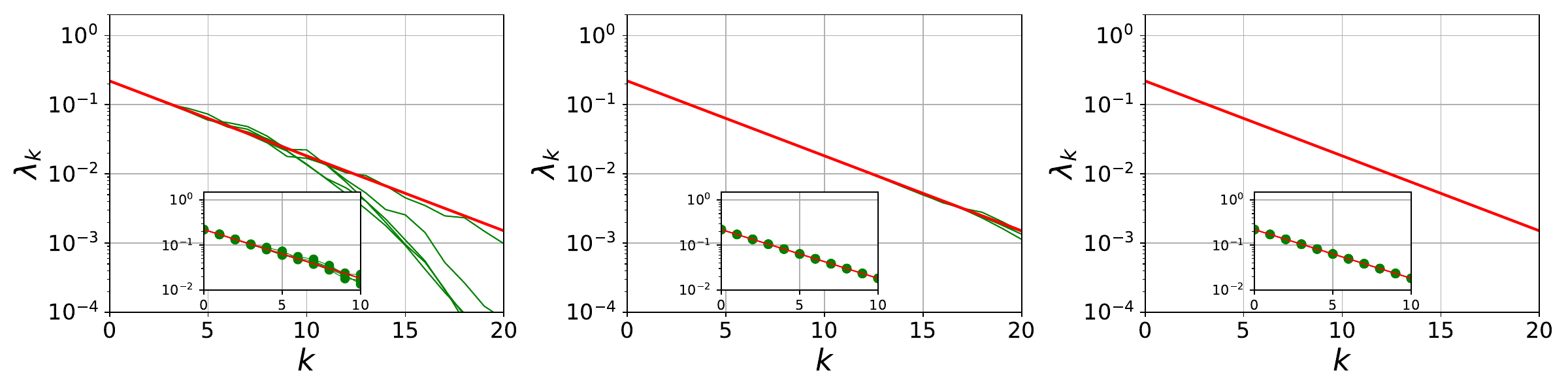}  
\includegraphics[width=0.85\linewidth]{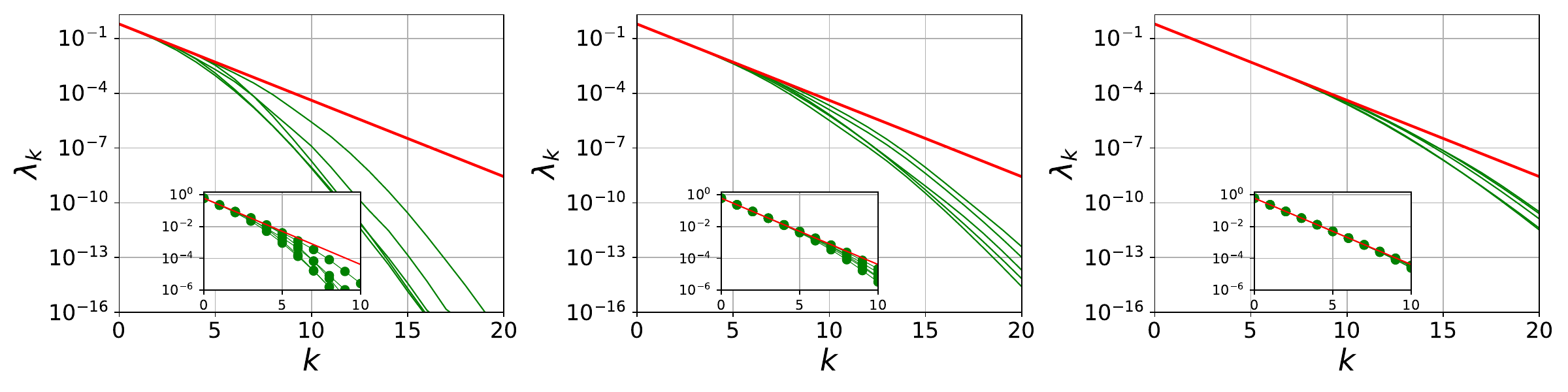}  
\end{center}
  \caption{\label{fig:ev_gauss_mc_s2} Karhunen-Lo\`{e}ve expansion eigenvalues for $\sigma_x=2$ and correlation lengths $\lc=\{\sigma_x/4,\sigma_x\}$ for top and bottom rows. The green lines show eigenvalues for different MC grid realizations. Other settings similar to Fig.~\ref{fig:ev_gauss_mc}.}
\end{figure}

\begin{figure}[!htb]
\begin{center}
\includegraphics[width=0.85\linewidth]{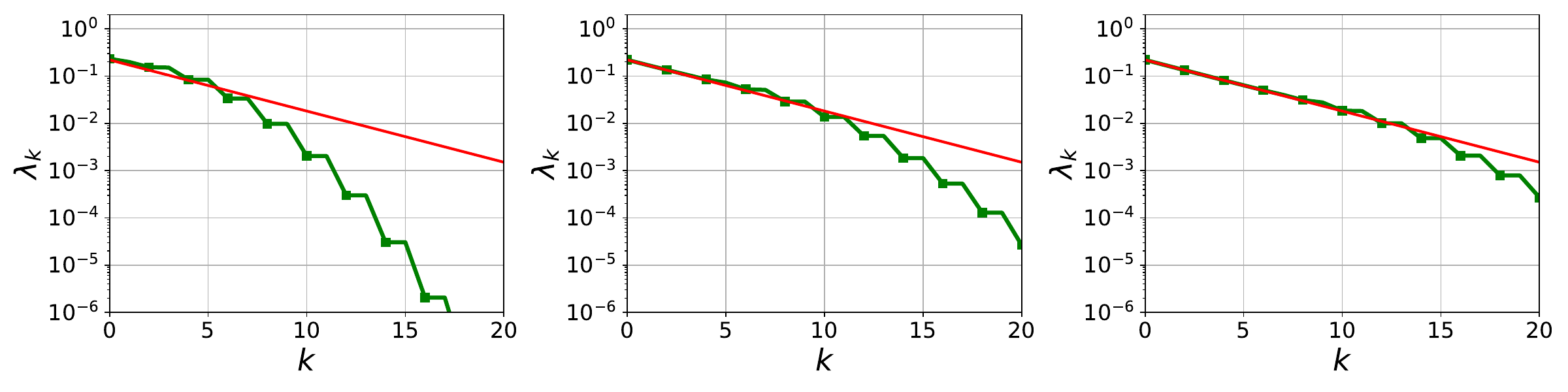}  
\includegraphics[width=0.85\linewidth]{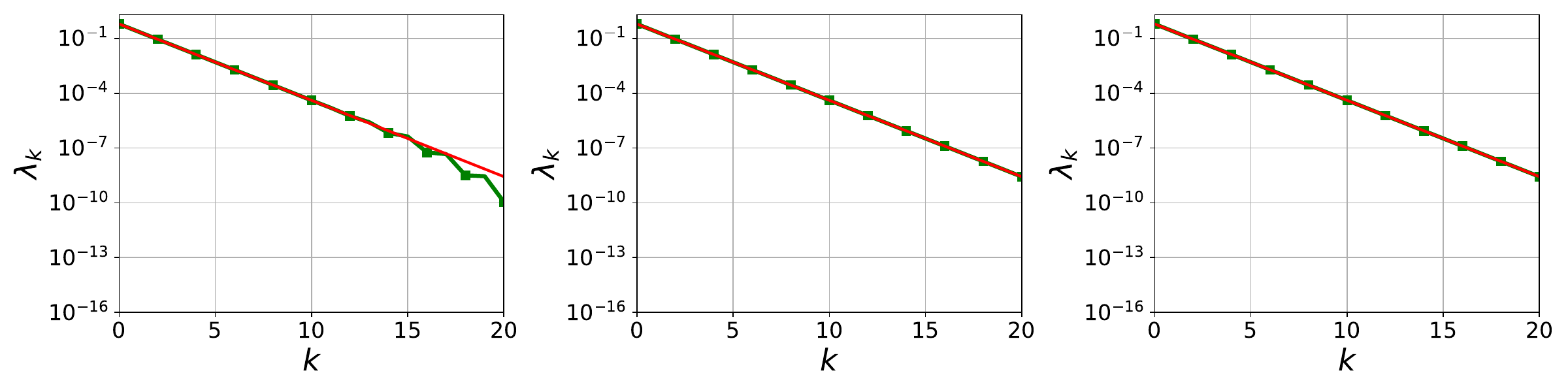}  
\end{center}
  \caption{\label{fig:ev_gauss_her_s2} Karhunen-Lo\`{e}ve expansion eigenvalues for a square exponential covariance matrix with the same set of correlation lengths as in Fig.~\ref{fig:ev_gauss_mc_s2}. Gauss-Hermite quadrature settings are the same as Fig.~\ref{fig:ev_gauss_her}.}
\end{figure}

To further verify the interplay between $\sigma_x$, $\lc$, and the number of quatrature points (MC or GH), Figs.~\ref{fig:ev_gauss_mc_s2} and~\ref{fig:ev_gauss_her_s2} show results for a case with $\sigma_x=2$, and $\lc/\sigma_x=\{1/4,1\}$, respectively. The MC-based results use different random seeds than for the results presented above and are qualitatively similar as the corresponding results for $\sigma_x=1$. The GH-based results are the same up to machine precision confirming the dependence on the ratio $\lc/\sigma_x$.  

\begin{figure}[h]
\begin{center}
\includegraphics[width=0.75\linewidth]{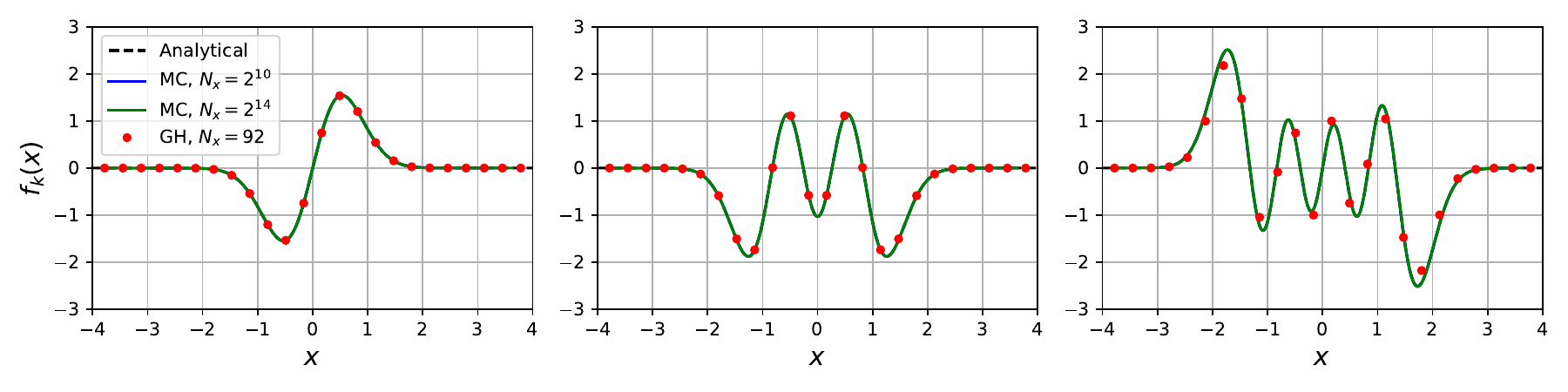}  
\includegraphics[width=0.75\linewidth]{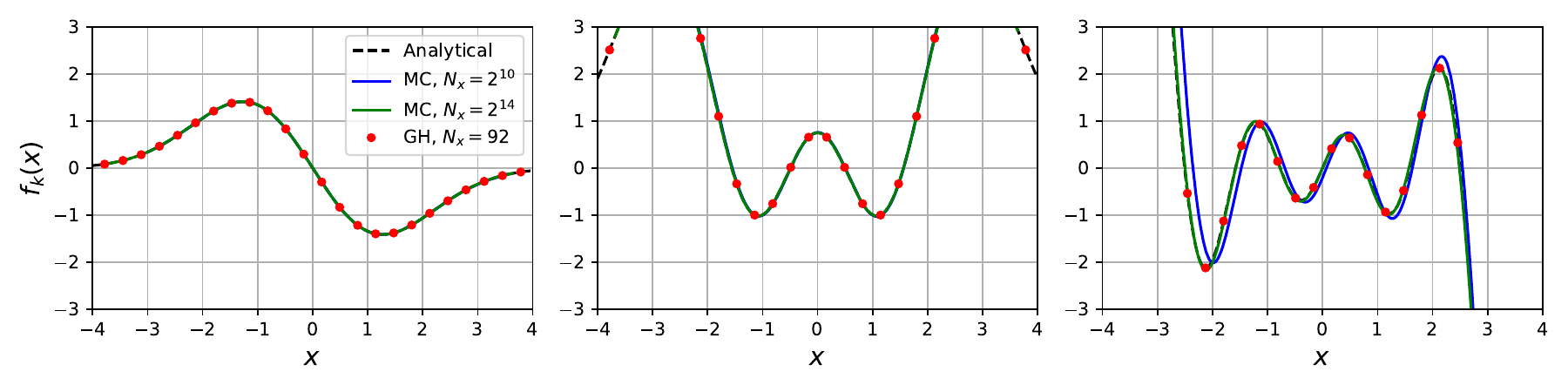}  
\includegraphics[width=0.75\linewidth]{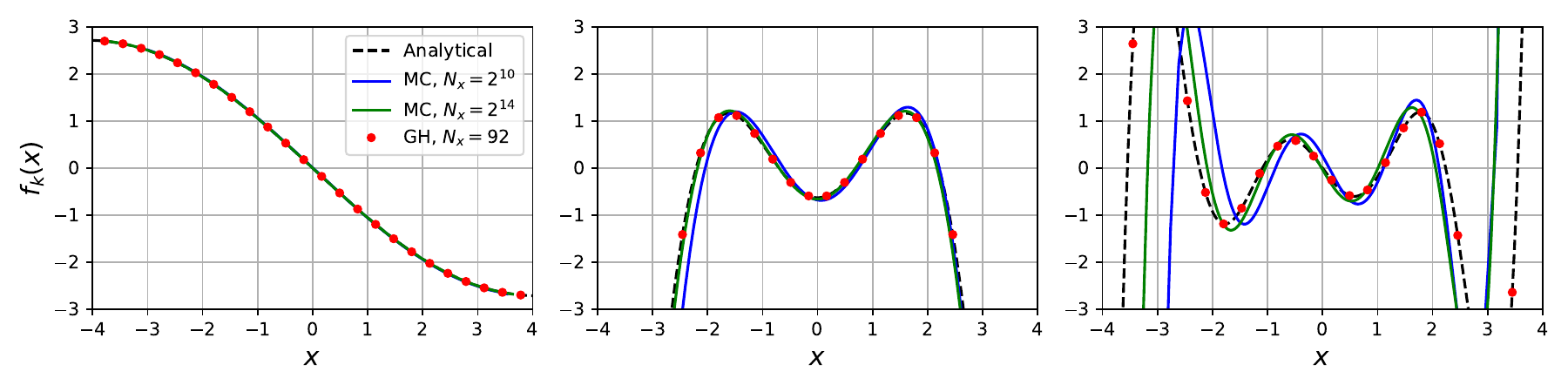}  
\end{center}
\caption{\label{fig:evec_gauss_her} Karhunen-Lo\`{e}ve expansion eigenmodes for a square exponential covariance matrix with the same set of correlation lengths as in Fig.~\ref{fig:ev_gauss_mc}. Dashed lines correspond to analytical eigenfunctions, while symbols correspond to numerical solutions via GH quadrature represented on the grid for the 92-order scheme. The left, center, and right columns correspond to modes 2, 5, and 8, respectively.}
\end{figure}

Wrapping up this section, we inspect select eigenvectors in Fig.~\ref{fig:evec_gauss_her}. The MC results for several grid sizes correspond to one of seeds shown in Fig.~\ref{fig:ev_gauss_mc}.  As $\lc$ increases, the modes become progressively wider and smoother, reflecting the underlying covariance structure.  Qualitatively, the agreement between numerical and analytical eigenfunctions presented in this figure correlates with interplay between $\lc$, $\sigma_x$, and $\Nx$. At small correlation lengths (top row), the MC-based results can outperform the GH results (albeit with a much larger number of grid samples - so the comparison is not one-to-one), given the concentration of MC in a narrow region. At the other end of the spectrum, for large $\lc$ in the bottom row, the GH quadrature grids extend well into the tail of the spatial distribution and reproduce well the eigenfunction shown in this plot. For this case, the MC results start showing visible discrepancies for higher order modes.

\clearpage

\subsection{Two-dimensional Karhunen-Lo\`{e}ve Expansion}

\begin{figure}[htb!]
\centering
\includegraphics[width=0.85\linewidth]{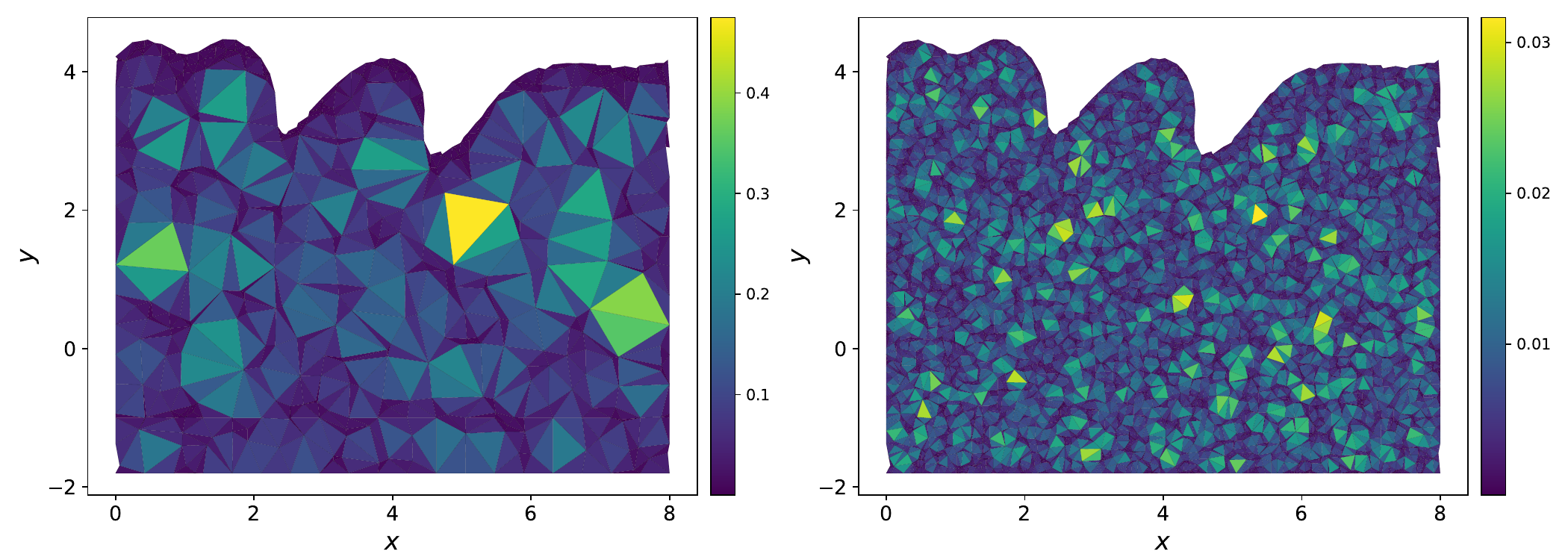}
\caption{\label{fig:2Dmesh_comp}Element areas (quadrature weights $w_j$) for the coarse mesh (left) and the dense mesh (right). The colorbars use different scales: the coarse mesh has areas up to $0.47$, while the dense mesh triangles are about $10\times$ smaller.}
\end{figure}

This section presents KLE results on a two-dimensional domain with a wavy top boundary. This domain is covered by an unstructured mesh with triangle elements. Figure~\ref{fig:2Dmesh_comp} shows the triangles for a coarse mesh in the left frame and a denser mesh in the right frame. Each triangle is colored by its area which also represents the quadrature weight of that element. The centers of these triangles represent the quadrature points that will be used to discretize the Fredholm equation. Four correlation lengths, $\lc\in[0.2,0.5,1,2]$ are used for this study.

The coarse mesh features approximately 500 vertices and 900 triangle elements, with a maximum triangle area of $A_{\max}=4.75\times10^{-1}$. The square root of this quantity, $\sqrt{A_{\max}}=0.689$ sets the effective resolution limit: correlation lengths smaller than this value cannot be adequately resolved by the coarse mesh. The fine mesh has about $10$  times more vertices and elements, $5000$ and $9500$, compared to the coarse mesh. The maximum element area is reduced by a factor of $15$ relative to the coarse mesh, and $\sqrt{A_{\max}} = 0.178$, which is below the smallest correlation length $\lc=0.2$ that will be tested in this section. 

Using the triangle centers as the quadrature points, we use Eq.~\eqref{eq:sqexpcov} to construct the covariance matrix and then employ the discretized Fredholm equation, presented in Sec.~\ref{sec:klenum}, to compute the eigenvalues and eigenvectors corresponding to a KL expansion on this irregular domain. The numerical quadrature is accurate provided that the mesh spacing is small relative to the correlation length.  For a triangular mesh a natural measure of local spacing is $\Delta_j = \sqrt{A_j}$. The midpoint rule can only represent the spatial changes of the random field within an element if there are sufficient elements to capture the fluctuations~\cite{Betz:2014,Ghanem:1991,Schwab:2006}, i.e. $\lc > \Delta_{\max} = \sqrt{A_{\max}}$.  

\begin{figure}[htb!]
\centering
\includegraphics[width=0.5\linewidth]{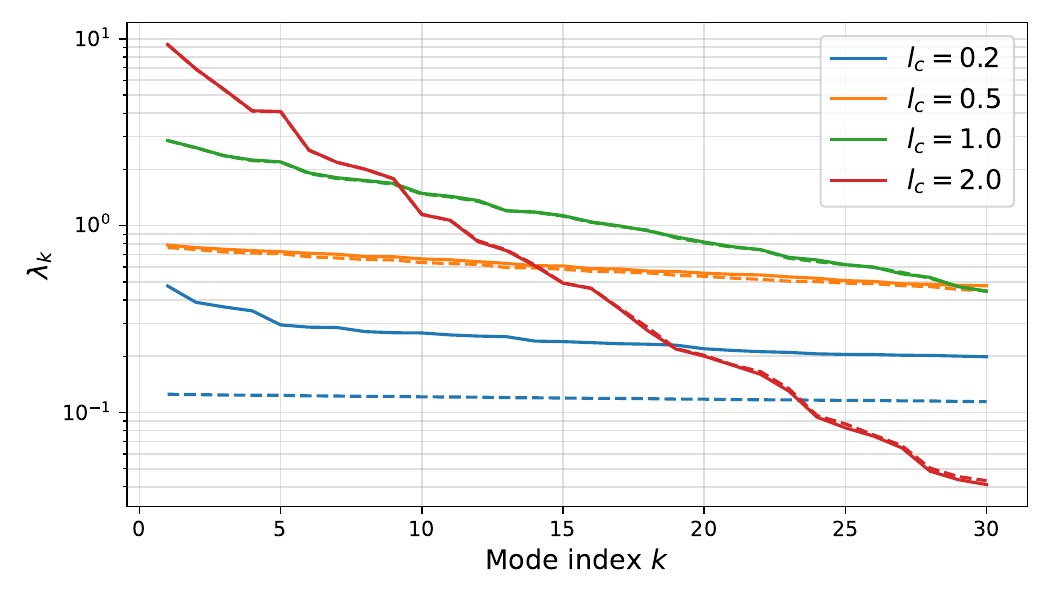}
\caption{\label{fig:evals_2D_unstr}KLE eigenvalue spectra for $\lc \in \{0.2,\,0.5,\,1.0,\,2.0\}$. Solid lines: coarse mesh (895 elements); dashed lines: dense mesh (9\,590 elements).}
\end{figure}

Figure~\ref{fig:evals_2D_unstr} shows the first 30 eigenvalues $\lambda_k$ on a semi-logarithmic scale. As expected, longer correlation lengths concentrate variance in fewer modes. At $\lc=2$, about 16 modes are needed to capture $95\%$ of the variance, while $\lc=1$ requires about 29 modes. 
For $\lc\in\{1,2\}$ the coarse and fine mesh spectra are virtually indistinguishable. At $\lc = 0.2$ the coarse mesh is under-resolved leading to eigenvalues significantly different than the fine mesh solution.

We now examine the eigenmodes for the cases with $\lc\geq 0.5$. To visualize the eigenmodes we employ barycentric interpolation to interpolate values from triangle centers to triangle vertices. Further, mode signs are aligned between the two meshes: for each mode $k$ the element with the largest $|\fs_k|$ in the coarse mesh is identified, the nearest element in the dense mesh is then located, and the dense-mesh eigenvector is flipped if the two signs disagree.

\begin{figure}[htbp]
\centering
\includegraphics[width=0.8\linewidth]{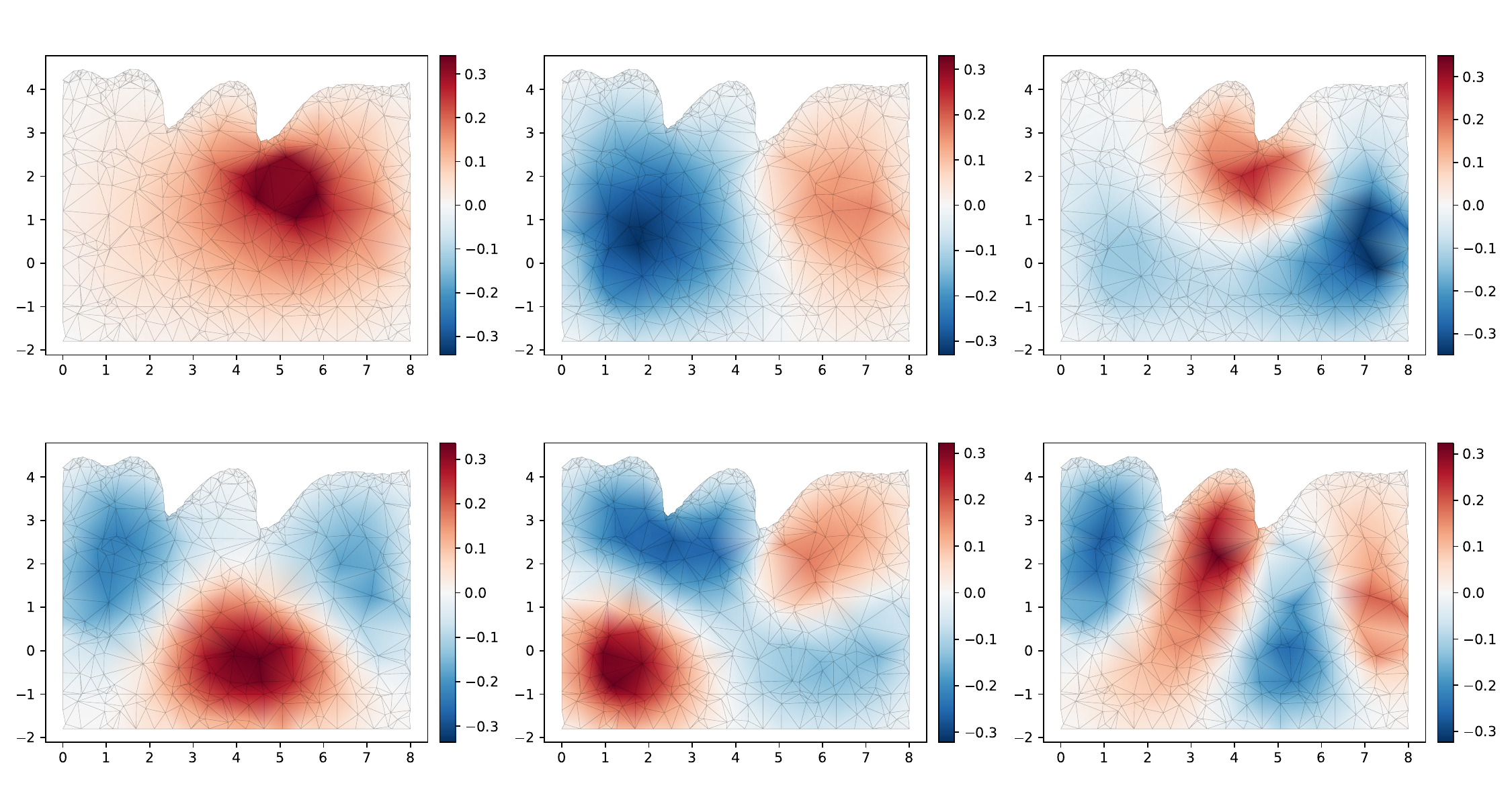}
\caption{\label{fig:modes_05_orig} Eigenmodes $k=1$-$6$ (left to right and top to bottom) for $\lc=0.5$, coarse mesh. The diverging colormap is centered at zero.}
\end{figure}

\begin{figure}[htbp]
\centering
\includegraphics[width=0.8\linewidth]{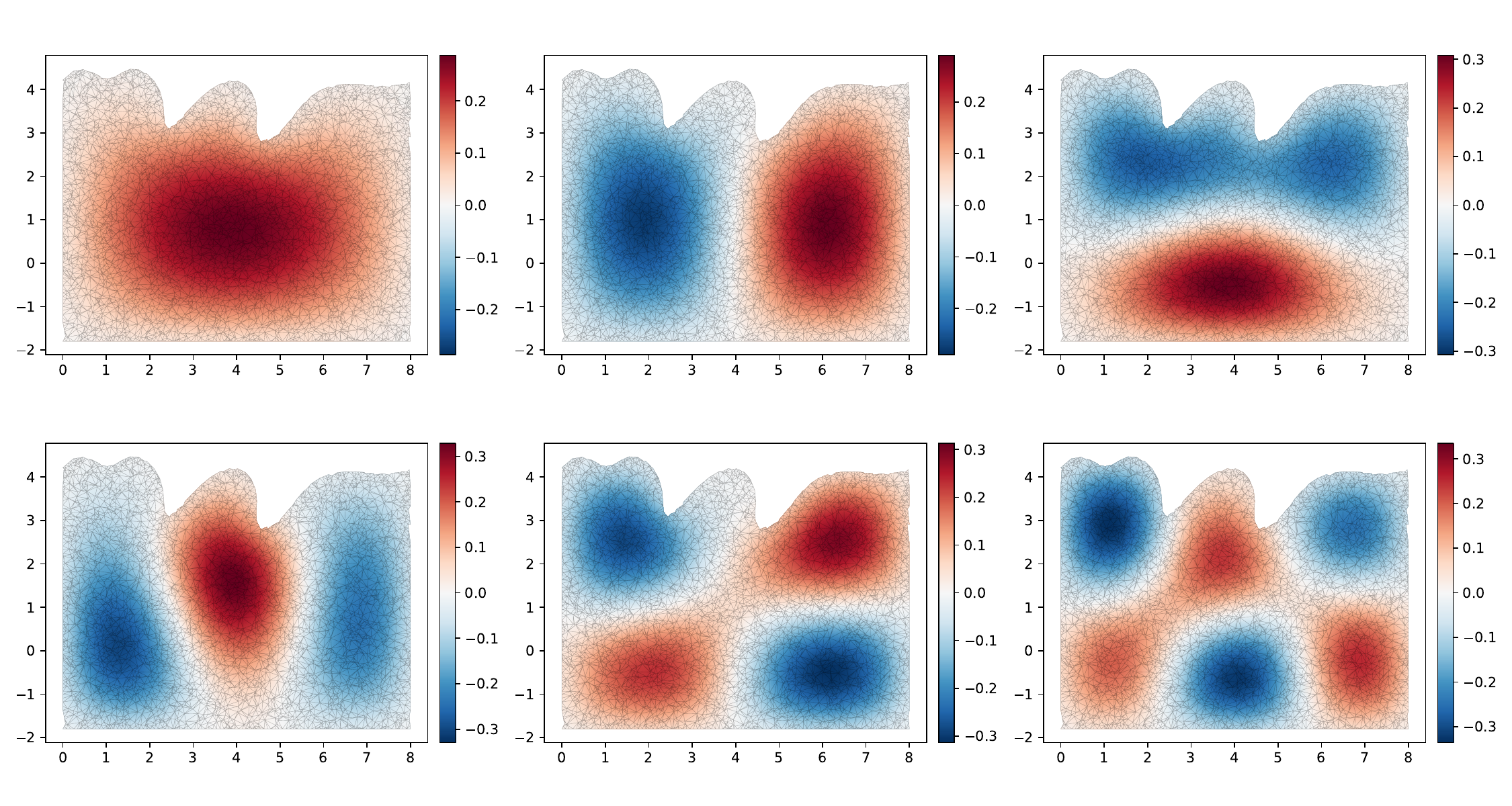}
\caption{\label{fig:modes_05_unif}Eigenmodes $k=1$-$6$ for $\lc=0.5$, dense mesh.
The spatial structures signs were matched to those in Figure~\ref{fig:modes_05_orig}.}
\end{figure}

Figures~\ref{fig:modes_05_orig} and~\ref{fig:modes_05_unif} show the first six modes for $\lc=0.5$ on the coarse and dense meshes respectively.
Mode~1 is a broad pole centered in the middle. Subsequent modes gradually introduce multipolar patterns. The spatial structures are qualitatively consistent between the two meshes.  Quantitative differences are observed, given the still large triangle areas for the coarse mesh relative to the correlation length for this case. In particular, mode 3 in the upper right frame could not be matched sign-wise between the two grids with the proposed algorithm given significant discrepancies between the two cases.

\begin{figure}[htbp]
\centering
\includegraphics[width=0.8\linewidth]{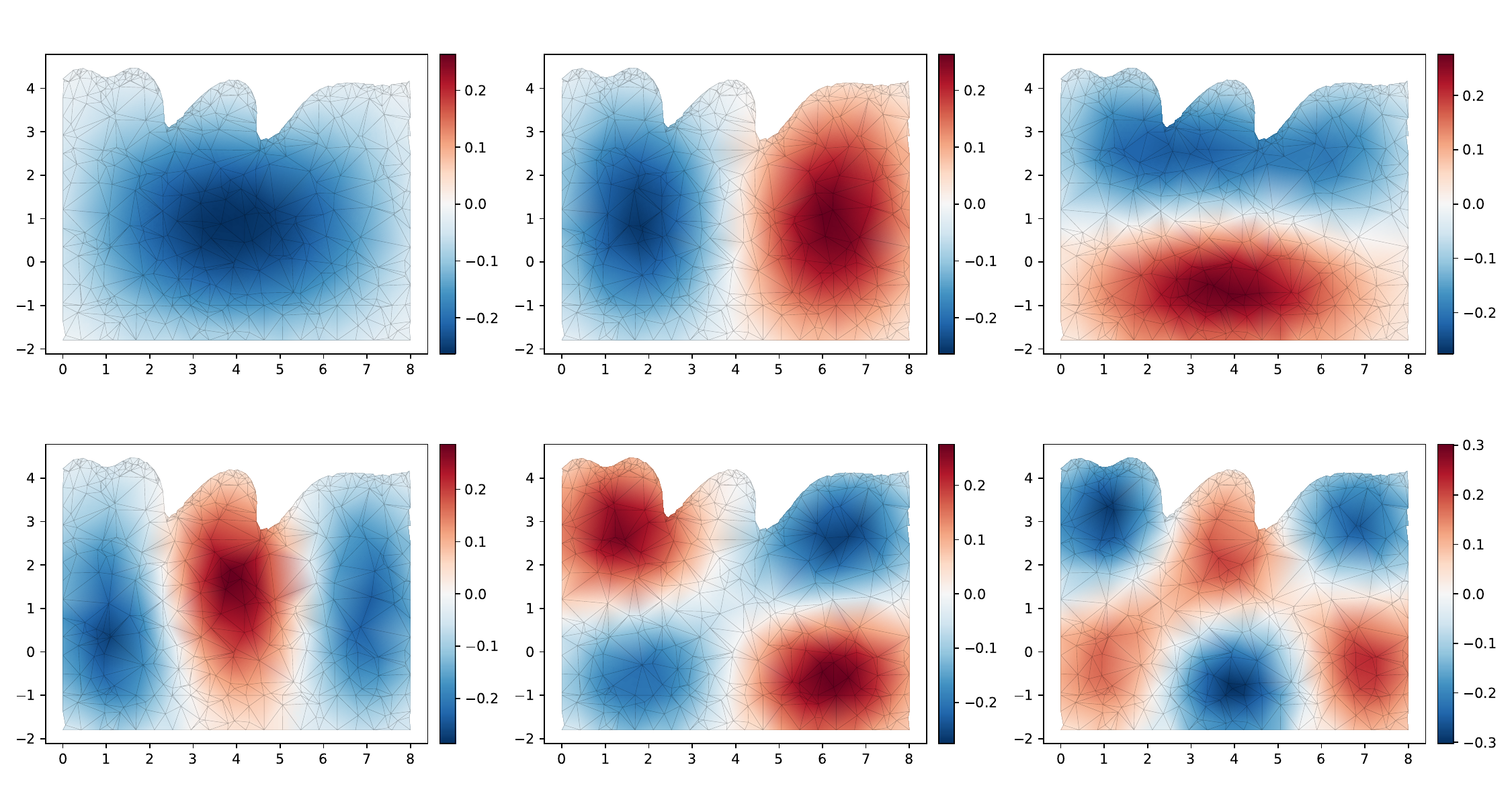}
\caption{\label{fig:2D_modes_10_orig}Eigenmodes $k=1$-$6$ for $\lc=1.0$, coarse mesh.}
\end{figure}

\begin{figure}[htbp]
\centering
\includegraphics[width=0.8\linewidth]{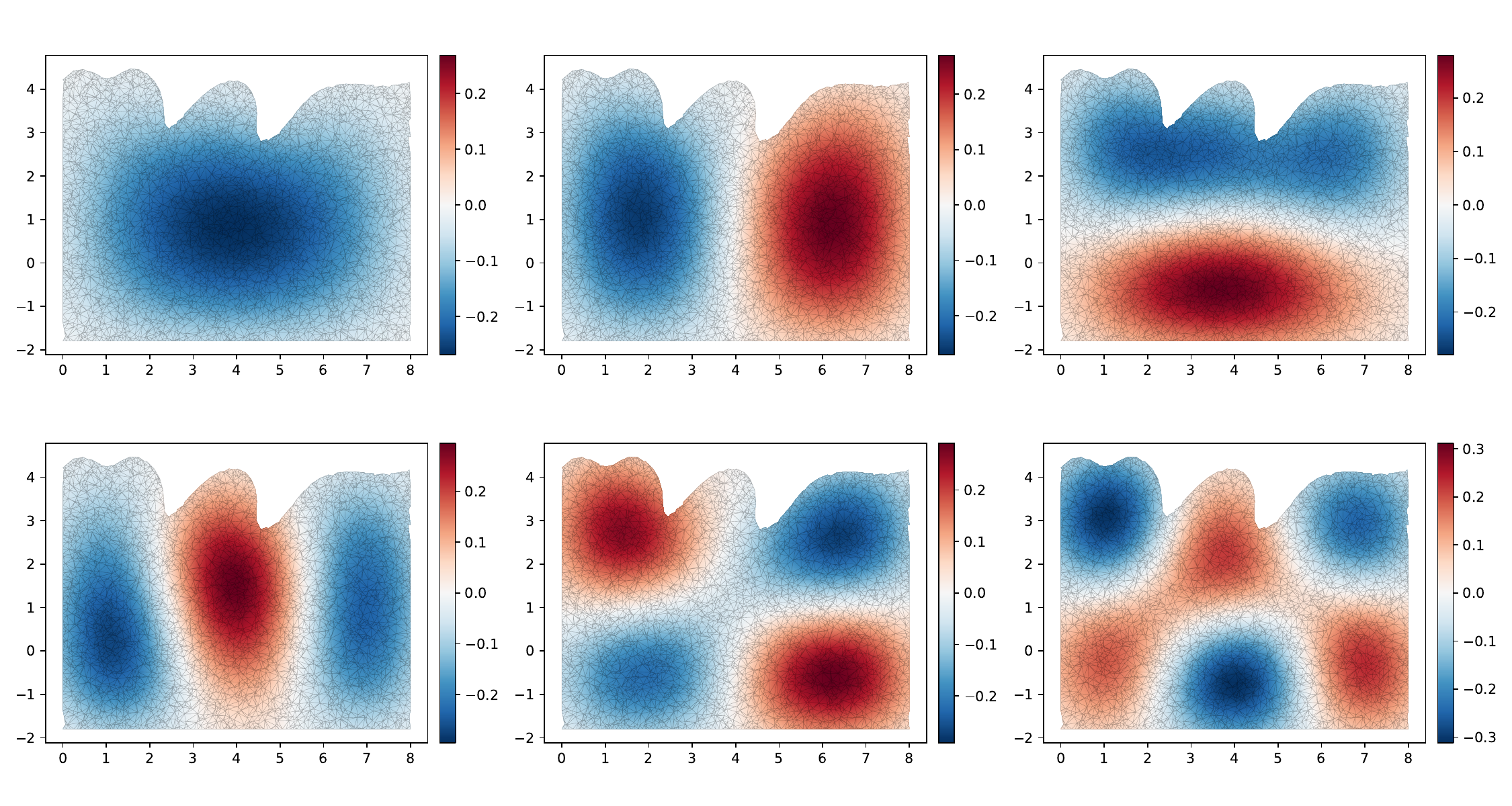}
\caption{\label{fig:2D_modes_10_unif}Eigenmodes $k=1$-$6$ for $\lc=1.0$, dense mesh. Structures are in close agreement with Figure~\ref{fig:2D_modes_10_orig}, consistent with the negligible eigenvalue discrepancies at this correlation length.}
\end{figure}

Figures~\ref{fig:2D_modes_10_orig} and~\ref{fig:2D_modes_10_unif} present a similar comparison for results corresponding to $\lc=1$. The two sets of grids produce visually identical structures, consistent with negligible errors present in the eigenvalue spectra. Results corresponding to $\lc=2$ (results not shown) are virtually indistinguishable.


\begin{figure}[htbp]
\centering
\includegraphics[width=0.8\linewidth]{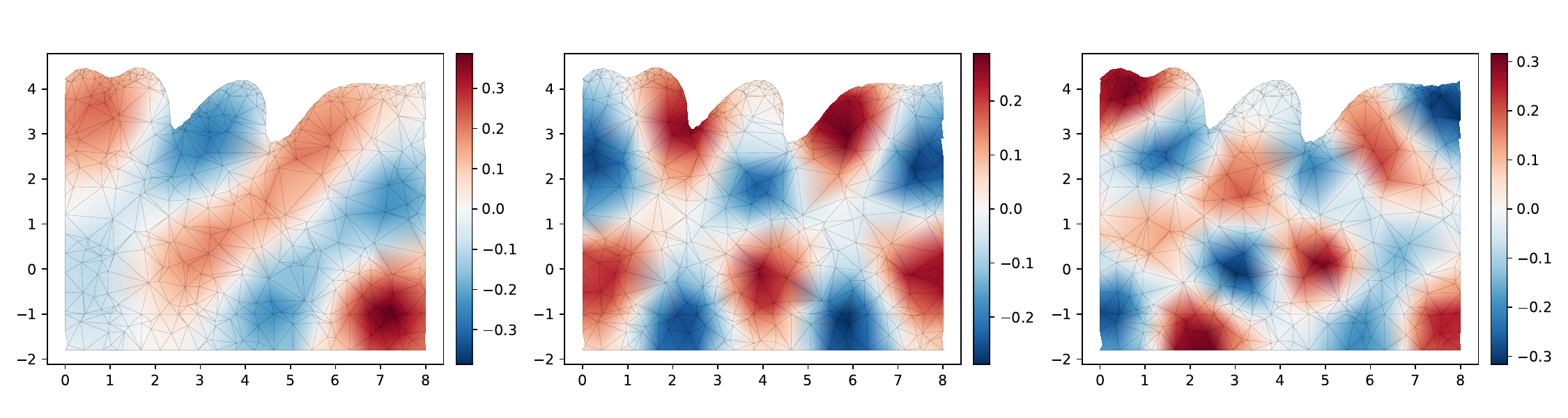}
\caption{\label{fig:2D_modes_hi_10}Modes $k=10$, $15$, $20$ for $\lc=1$, coarse mesh.}
\end{figure}

\begin{figure}[htbp]
\centering
\includegraphics[width=0.8\linewidth]{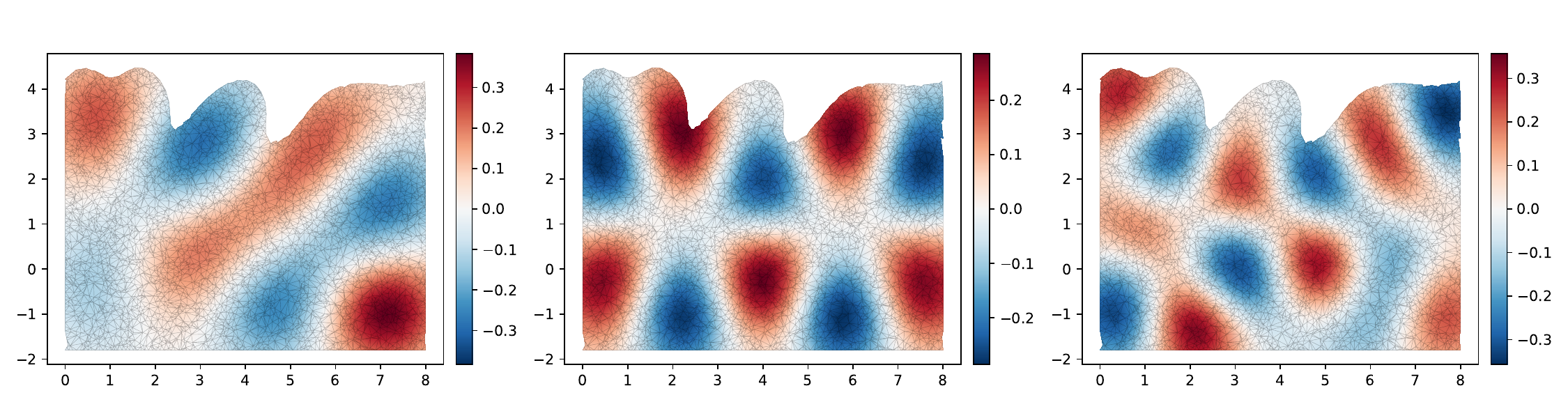}
\caption{\label{fig:2D_modes_hi_10_unif}Modes $k=10$, $15$, $20$ for $\lc=1$, dense mesh.}
\end{figure}

Figures~\ref{fig:2D_modes_hi_10} and~\ref{fig:2D_modes_hi_10_unif} show modes 10, 15, and 20 for $\lc=1$. These modes exhibit a multipolar pattern with spatial oscillations becoming increasingly fine with mode index, as expected. The eigenfunction amplitudes remain broadly similar across these three modes.  The dense mesh renders all three modes more smoothly but the overall spatial topology is reasonably reproduced by the coarse mesh, confirming adequate resolution for $\lc=1.0$.

\subsection{Three-dimensional Karhunen-Lo\`{e}ve Expansion}

In this section we consider random fields defined over the interior of a torus
\begin{equation}
\mathfrak{D} = \bigl\{\mathbf{x}\in\mathbb{R}^3 \;:\;
\bigl(\sqrt{x^2+y^2}-R\bigr)^2 + z^2 \leq r^2 \bigr\},
\end{equation}
with major radius $R=3$ (distance from the torus axis to the tube center) and minor (tube) radius $r=1$. The torus axis is the same as the z-axis for the coordinate sytem. We create a computational mesh for the torus as follows. First, a uniform $41\times41\times41$ Cartesian grid is constructed over the bounding box $[-4.05,4.05]^2 \times [-1.05,1.05]$, resulting in computational mesh sizes $\Delta_x = \Delta_y = 0.2025$, $\Delta_z=0.0525$. Only cells with one of the eight corners satisfying $(\sqrt{x^2+y^2}-R)^2+z^2\leq r^2$ are used in the computation.  After thresholding, $N\approx 22\times 10^3$ cells remain inside the torus.

\begin{figure}[htb!]
\centering
\includegraphics[width=0.6\linewidth]{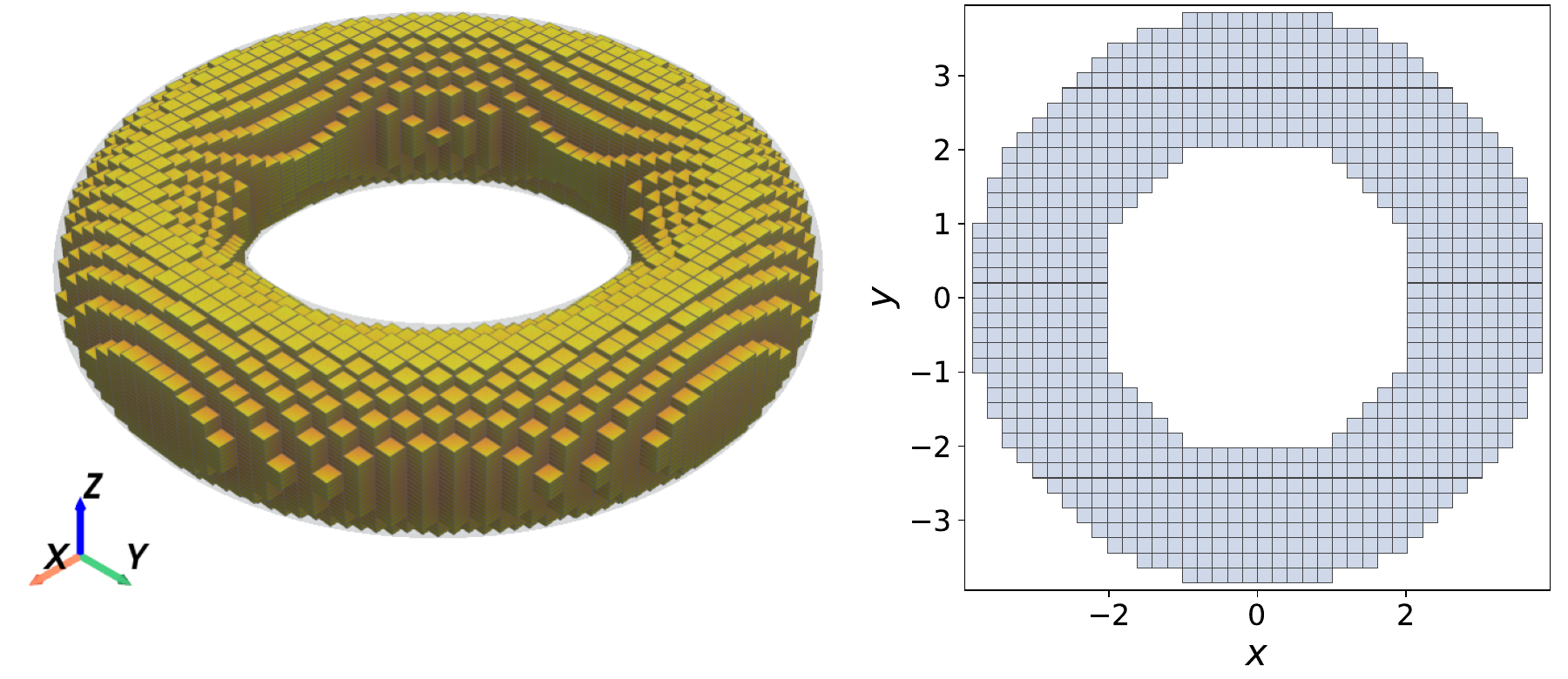}
\caption{\label{fig:3D_geometry}Structured Cartesian grid clipped to the torus interior. \textit{Left:} isometric 3-D view of the visible cell blocks. \textit{Right:} $z=0$ cross-section showing the individual Cartesian cells as rectangles, forming an annulus with inner radius $R-r=2$ and outer radius $R+r=4$.}
\end{figure}

Figure~\ref{fig:3D_geometry} shows the Cartesian grid cells clipped to the torus interior. The right panel shows the $z=0$ horizontal cross-section: each cell that intersects the mid-plane is rendered as an individual rectangle, forming the annular ring with inner radius $R-r=2$ and outer radius $R+r=4$.

Standard covariance kernels depend on the \emph{Euclidean} distance between two points.  When the domain is not simply connected - as is the case for the torus domain in this section - a pair of points separated by the central hole may have a short Euclidean distance but a long path distance \emph{through the material}. The \emph{shortest interior path} (SIP) distance, defined as the length of the shortest path confined to the domain interior, is then the physically more appropriate argument for the covariance function as the correlations between points separated by the hole are correctly attenuated by the true distance through the material.

\begin{figure}[htb!]
\centering
\includegraphics[width=0.7\linewidth]{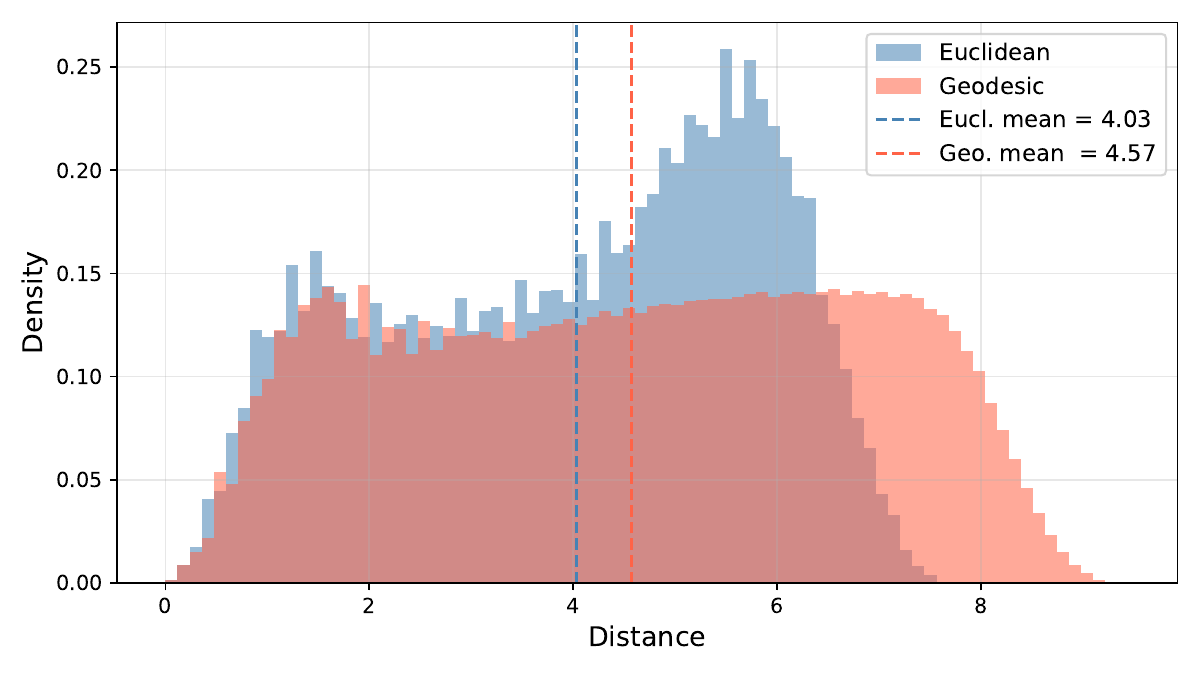}
\caption{\label{fig:3D_dists}Empirical density of pairwise distances between the $N=22\,168$ cell centroids, estimated from $2\times10^6$ randomly sampled pairs.  Blue: Euclidean distances; red: SIP distances.  Dashed vertical lines mark the respective means. The later distribution is right-shifted because pairs on opposite sides of the central hole are connected around the hole.}
\end{figure}

SIP distances are numerically estimated by constructing a sparse nearest-neighbor graph on the $N$ cell centroids and solving the all-pairs shortest-path problem with Dijkstra's algorithm.  All-pairs shortest paths are computed via \texttt{scipy.sparse.csgraph.shortest\_path}, yielding the symmetric $N\times N$ distance matrix $\mathsf{D}$. Figure~\ref{fig:3D_dists} overlays the empirical distributions of pairwise Euclidean and SIP distances.
The Euclidean distribution has a peak near $d\approx 5.5$, reflecting pairs roughly across the torus ring.  The SIP distribution has a longer right tail: for pairs separated by the hole, the interior path must travel longer around the ring. The mean SIP distance exceeds the mean Euclidean distance by approximately $10\%$, with the excess concentrated in the tail corresponding to cross-hole pairs.

In this section we will present results corresponding to a squared-exponential kernel, with distances between points being either the Euclidean distance or the SIP distance. For both cases we will consider three correlation lengths $\lc\in\{0.5,1,2\}$.

\begin{figure}[htb!]
\centering
\includegraphics[width=0.6\linewidth]{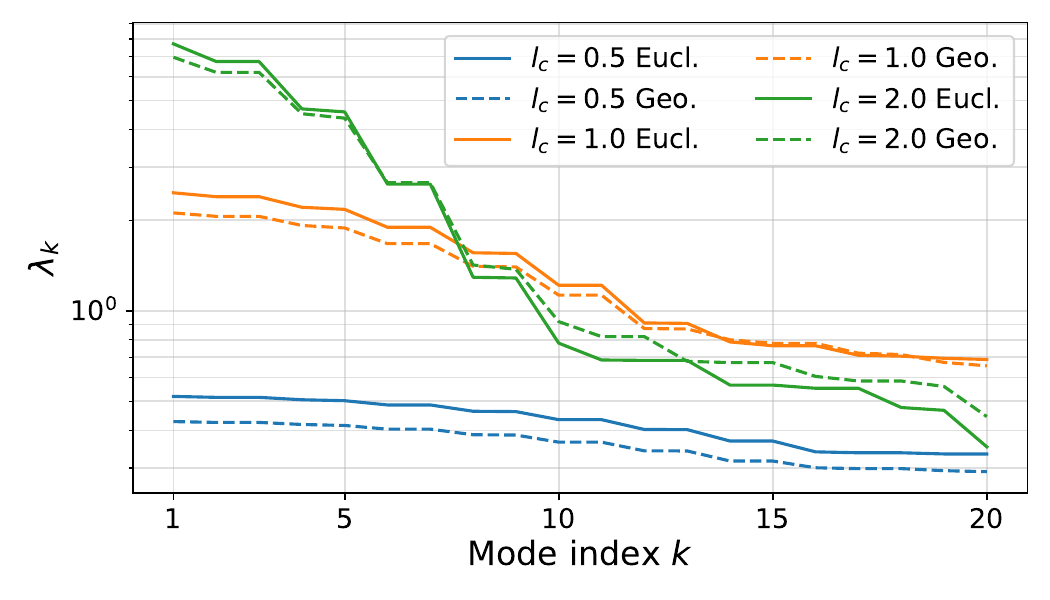}
\caption{\label{fig:evals_3D}KL expansion eigenvalue spectra for $\lc\in\{0.5,\,1.0,\,2.0\}$. Solid lines: Euclidean kernel; dashed lines: SIP kernel.}
\end{figure}

Figure~\ref{fig:evals_3D} shows the first 20 eigenvalues $\lambda_k$ on a
semi-logarithmic scale for both kernels and all three correlation lengths.
The effects of correlation length magnitude on the eigenvalue decay is similar to other cases presented above. For every $\lc$ the SIP leading eigenvalue is smaller than the Euclidean one and the SIP kernel attenuates the variance decay, redistributing variance across higher modes. 

\begin{figure}[htb!] 
\centering
\includegraphics[width=0.8\linewidth]{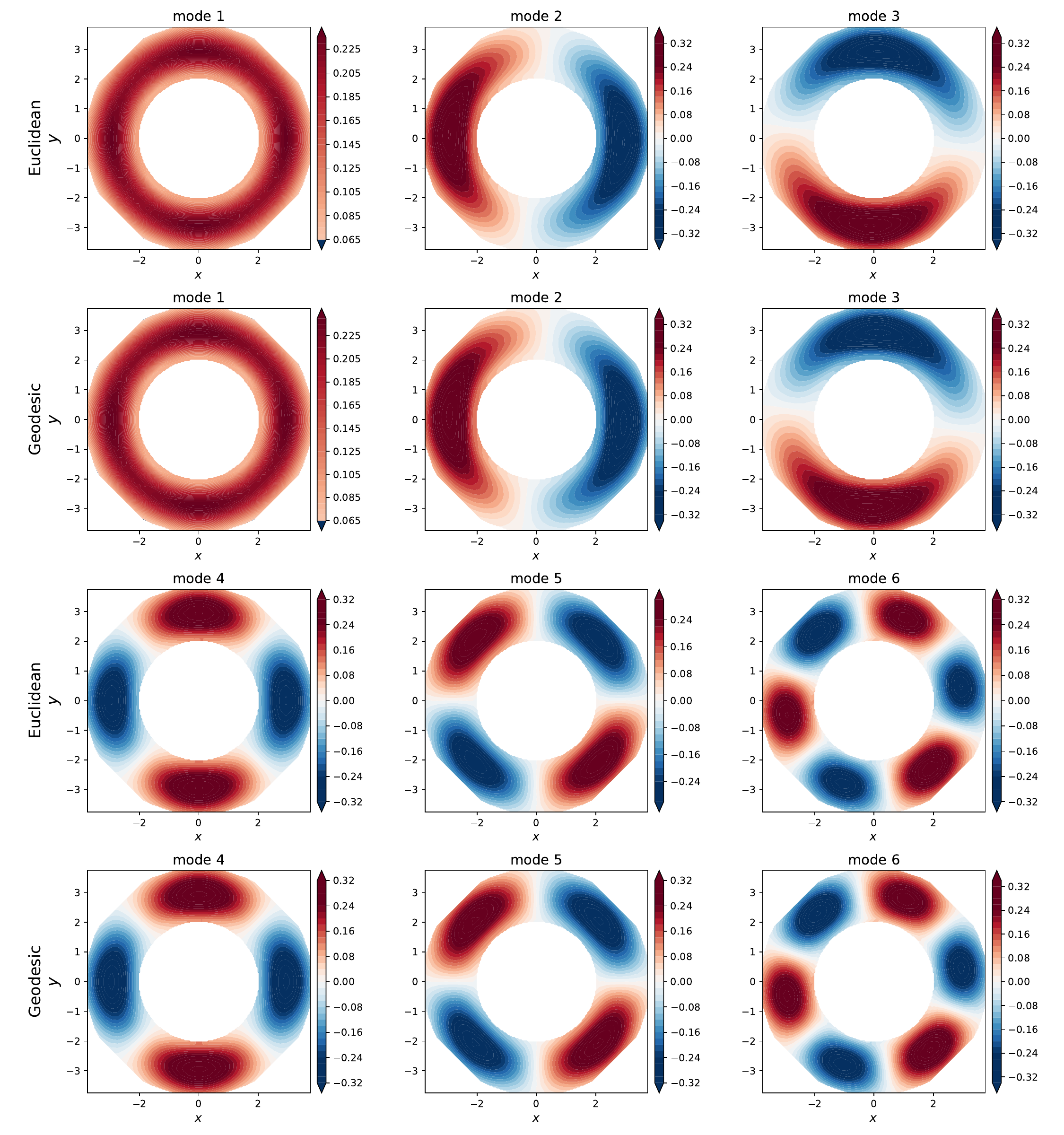}
\caption{\label{fig:3D_modes_z0_05}First six eigenmodes at $z=0$ for $\lc=0.5$. Rows 1 \& 3: Euclidean kernel.  Rows 2 \& 4: SIP kernel.}
\end{figure}

\begin{figure}[htb!]
\centering
\includegraphics[width=0.8\linewidth]{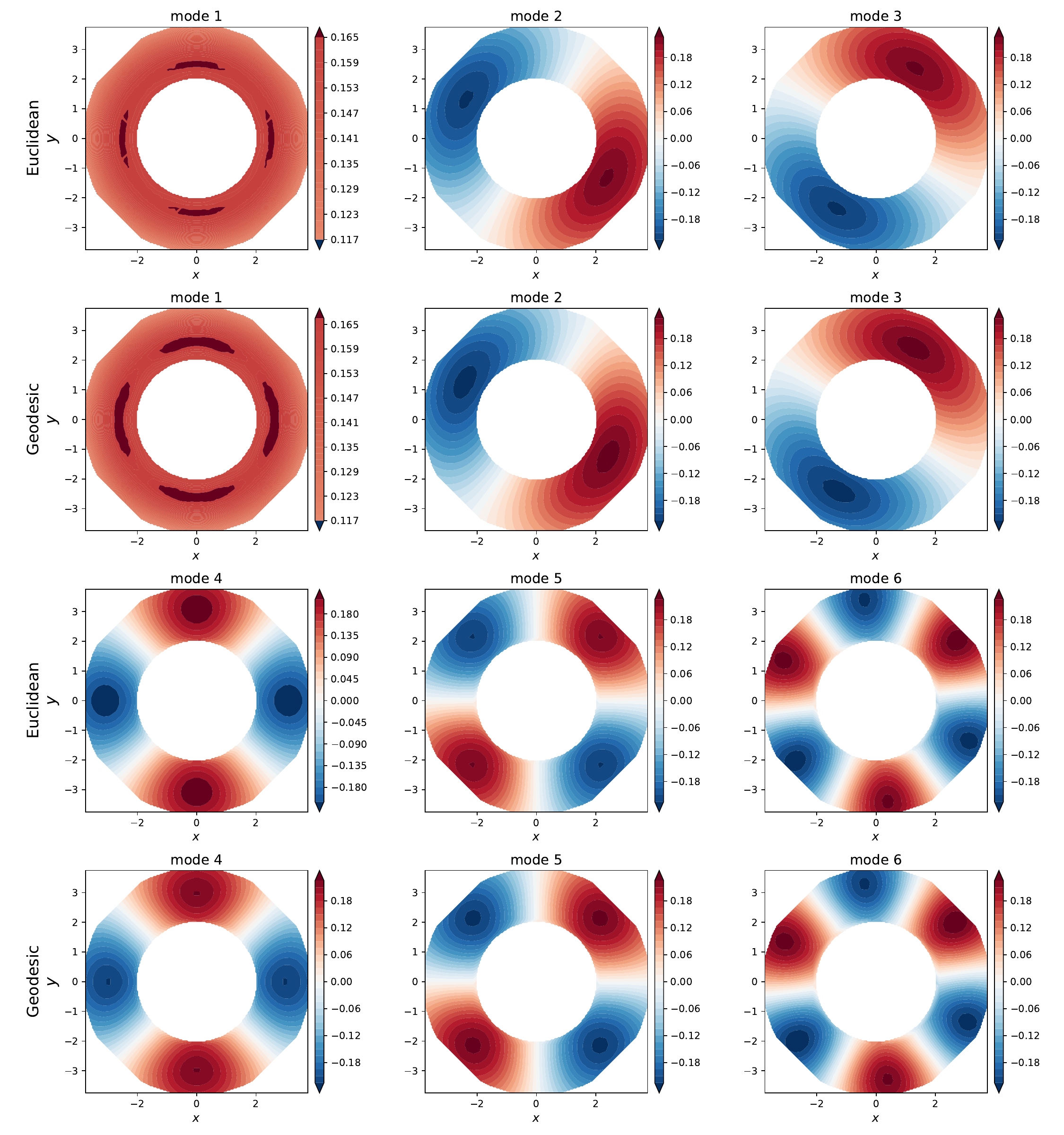}
\caption{\label{fig:3D_modes_z0_20}First six eigenmodes at $z=0$ for $\lc=2$. Rows 1 \& 3: Euclidean kernel.  Rows 2 \& 4: SIP kernel.}
\end{figure}

\begin{figure}[htb!]
\centering
\includegraphics[width=0.8\linewidth]{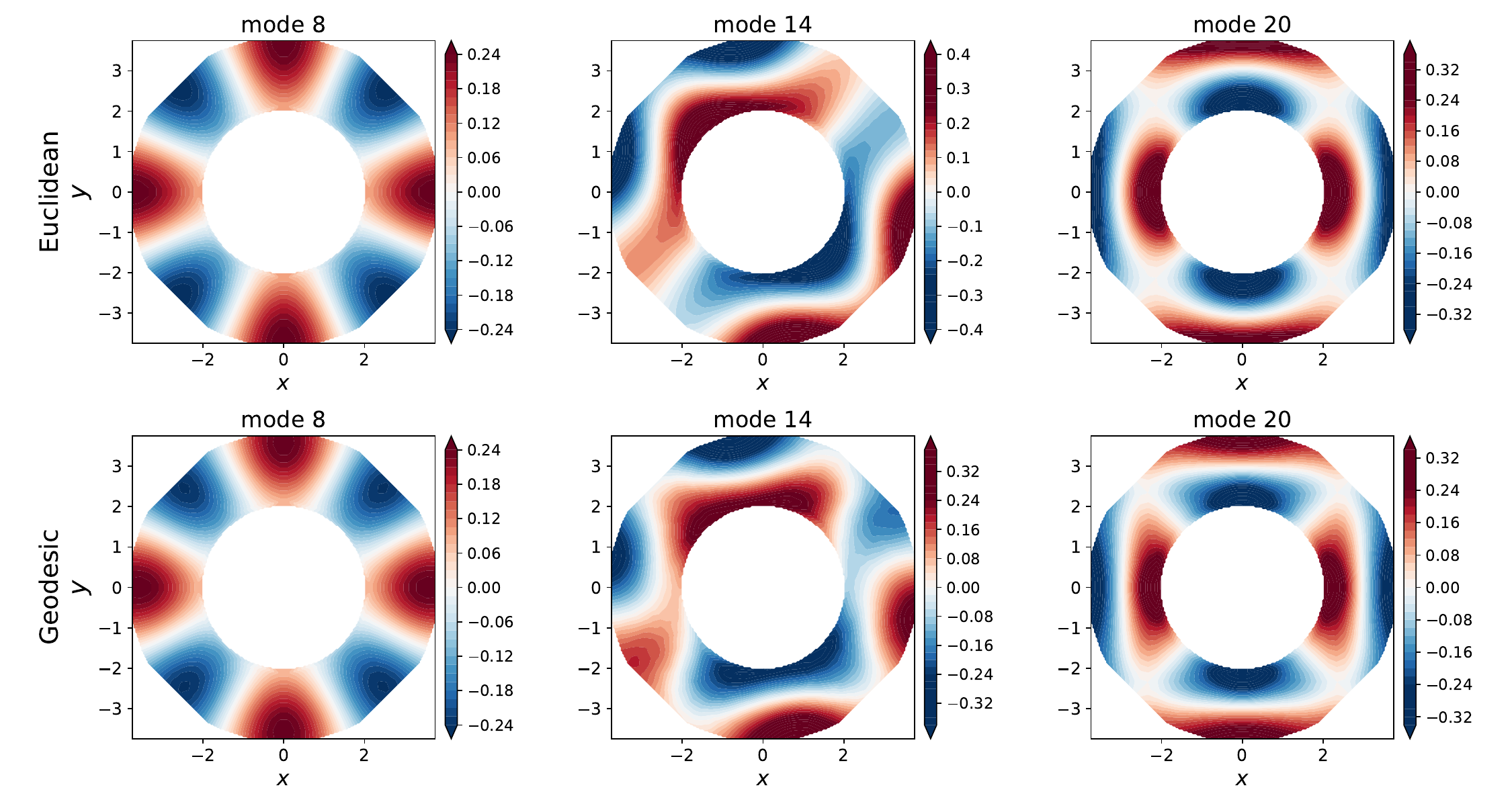}
\caption{\label{fig:3D_modes_z0_20_hi}Select higher eigenmodes at $z=0$ for $\lc=2$.}
\end{figure}

For visualization, eigenvectors are evaluated at cell centroids and interpolated onto a $200\times200$ regular grid using bilinear interpolation, in the $z=0$ plane, in Figures~\ref{fig:3D_modes_z0_05}, \ref{fig:3D_modes_z0_20}, and \ref{fig:3D_modes_z0_20_hi}.  Only the cells inside the torus cross-section are rendered in the plots below. Each figure combines Euclidean and SIP modes as follows: rows 1 and 2 correspond to modes 1-3 and rows 3 and 4 show modes 4-6. The modes are sign- and angle-aligned to enable a direct comparison between the two sets of results.

The $z=0$ section intersects the full tube cross-section. Figures~\ref{fig:3D_modes_z0_05} and~\ref{fig:3D_modes_z0_20} show Modes 1-6 for $\lc=0.5$ and $\lc=2$, respectively.  Figure~\ref{fig:3D_modes_z0_20_hi} shows higher order modes for $\lc=2$.
The Euclidean and SIP modes are visually consistent in shape for all tested $\lc$. The most pronounced differences appear in modes with spatial support on opposite sides of the annulus simultaneously, as can be seen for Mode 20 in Figure~\ref{fig:3D_modes_z0_20_hi}, where the two kernels differ most in their cross-hole correlation.

\section{Summary}
\label{sec:conclusions}
This report examines numerical aspects of constructing Karhunen-Lo\`{e}ve expansions (KLEs) for second-order stochastic processes, with emphasis on the interplay between spatial discretization strategy, covariance kernel structure, correlation length, and the number of available random-field samples. Two complementary solution pathways were considered: the direct numerical solution of the discretized Fredholm integral equation of the second kind via a symmetrized eigenproblem and the algebraically equivalent singular value decomposition of the weight-scaled sample matrix.

In one dimension, numerical results for the exponential kernel confirm that a grid of given spacing accurately reproduces the analytical spectrum up to mode $k$ provided there are enough samples to capture the mode oscillations; eigenvalues are significantly more sensitive to under-resolution than eigenfunctions. The inter-seed variability of SVD eigenvalue estimates and the Kullback-Leibler divergence of the empirical KLE coefficient distributions from $\mathcal{N}(0,1)$, both converge at the Monte Carlo rate $\mathcal{O}(1/\sqrt{\Ns})$, uniformly across mode indices. For the squared-exponential kernel on $\mathcal{N}(0,\sigma_x^2)$, eigenvalues decay geometrically with a per-mode rate determined entirely by the ratio $\rho = \lc/\sigma_x$ between the correlation length and the standard deviation corresponding to the spatial density. Gauss-Hermite quadrature is superior for $\rho \gtrsim 1$ due to its extended tail coverage, while Monte Carlo sampling is competitive for $\rho<1$ where the dominant eigenfunctions are spatially compact.

On a two-dimensional unstructured triangular mesh, the resolution criterion $\lc > \Delta_\mathrm{max} = \sqrt{A_\mathrm{max}}$ is necessary for converged numerical solutions. The three-dimensional example features a multiply-connected domain, resulting in the relevance of chosen distance measures. The shortest interior path distance kernel is a physically consistent approach to estimate correlations between sample points. Differences between the shortest interior path and Euclidean kernels are negligible for low-order modes but become visible for higher mode indices.

\section*{Acknowledgement}
This work was supported in part by the U.S. Department of Energy, Office of Science, Office of Advanced Scientific Computing Research, Scientific Discovery through Advanced Computing (SciDAC) program through the FASTMath Institute. This report has been co-authored by employees of National Technology and Engineering Solutions of Sandia, LLC under Contract No. DE-NA0003525 with the U.S. Department of Energy (DOE).

\bibliographystyle{plain}

\end{document}